\documentclass[a4paper,
               10pt,
               pdftex,
               normalheadings,
               headsepline,
               footsepline,
               headinclude,
               footinclude,
               DIV=14,
               abstracton]
{scrartcl}

\usepackage[left=3cm, right=3cm, top=2.50cm, bottom=3.50cm]{geometry}

\usepackage
{
    graphicx,
    amssymb,
    amsmath,
    amsthm,
    xcolor,
    dsfont,
    algpseudocode,
    authblk,
    bm,
    bbm
}

\usepackage[ngerman, english]{babel}
\usepackage{enumerate}
\usepackage{algorithm}
\usepackage{tabularx}
\usepackage{caption}
\usepackage{subcaption} % subfigures
\usepackage{booktabs} % subfigures
\usepackage{siunitx}
\usepackage{comment}

\usepackage[colorlinks,
            pdffitwindow=false,
            plainpages=false,
            pdfpagelabels=true,
            pdfpagemode=UseOutlines,
            pdfpagelayout=SinglePage,
            bookmarks=false,
            colorlinks=true,
            hyperfootnotes=false,
            linkcolor=blue,
            urlcolor=blue!50!black,
            citecolor=green!50!black]
{hyperref}
\usepackage{cleveref}

\usepackage{sectsty}
\sectionfont{\fontsize{11}{11}\selectfont}
\subsectionfont{\fontsize{10}{10}\selectfont}

\newtheorem{theorem}{Theorem}[section]
\newtheorem{corollary}[theorem]{Corollary}
\newtheorem{lemma}[theorem]{Lemma}
\newtheorem{proposition}[theorem]{Proposition}
\newtheorem{definition}[theorem]{Definition}
\newtheorem{remark}[theorem]{Remark}
\newtheorem{example}[theorem]{Example}

\newtheorem{conjecture}[theorem]{Conjecture}

\let\OLDthebibliography\thebibliography
\renewcommand\thebibliography[1]{
	\OLDthebibliography{#1}
	\setlength{\parskip}{0pt}
	\setlength{\itemsep}{0pt plus 0.3ex}
}

\newcommand{\R}{\mathbb{R}}
\newcommand{\Pb}{\mathbb{P}}
\newcommand{\rv}[1]{\bm{#1}}
\newcommand{\super}[1]{^{{#1}}}
\newcommand{\set}[1]{\mathcal{#1}}
\newcommand{\toprob}{\overset{\mathrm{p}}{\longrightarrow}}
\newcommand{\diff}{\mathrm{d}}
\newcommand{\ldif}{l}

\newcommand{\review}[1]{{#1}}

\begin{document}

\title{Large population limits of Markov processes on random networks}
\author[1]{Marvin Lücke}
\author[3]{Jobst Heitzig}
\author[2]{Péter Koltai}
\author[4]{Nora Molkenthin}
\author[1]{Stefanie Winkelmann}

\affil[1]{Zuse Institute Berlin}
\affil[2]{Department of Mathematics, Free University of Berlin}
\affil[3]{FutureLab on Game Theory and Networks of Interacting Agents, Potsdam Institute for Climate Impact Research}
\affil[4]{Complexity Science Department, Potsdam Institute for Climate Impact Research}

\date{}

\maketitle
\begin{abstract}
We consider time-continuous Markovian discrete-state dynamics on random networks of interacting agents and study the large population limit. The dynamics are projected onto low-dimensional collective variables given by the shares of each discrete state in the system, or in certain subsystems, and general conditions for the convergence of the collective variable dynamics to a mean-field ordinary differential equation are proved. We discuss the convergence to this mean-field limit for a continuous-time noisy version of the so-called ``voter model'' on Erd\H{o}s--Rényi random graphs, on the stochastic block model, and on random regular graphs. Moreover, a heterogeneous population of agents is studied.
\end{abstract}

\section{Introduction}
Large networks, where the nodes represent agents and the edges indicate some form of interaction between them, are used to model numerous (social) phenomena \cite{Easley2010, Castellano2009}, e.g., the spreading of a disease \cite{Kiss2017} or the diffusion of a certain (political) opinion within a society \cite{das2014}, just to name a few.
In such a framework each node has a state that changes over time, in dependence upon the states of other nodes.
Stochastic effects are often included in the model in order to take account of uncertainty in the dynamics and variability of the agents' behavior \cite{Castellano2009}.  Even for simple interaction rules dictating the evolution of an agent's state, the emergent macroscopic system behavior remains difficult to examine analytically. Hence, one usually resorts to numerical simulations in order to analyze such systems \cite{Porter2016}.
As the computational effort for running simulations typically increases linearly or faster with the size of the network, simulating large populations of agents is often cumbersome or even infeasible. 
Thus, while on the one hand modeling each agent’s behavior simultaneously carries the hope of capturing effects that are elusive to less detailed models, on the other hand one is interested in reduced-order models that retain the central phenomena of the detailed model and still allow for a good understanding and computational feasibility.
One approach to address this issue is to find a low-dimensional representation of the system which captures the dynamics in the large population limit.

One of the classical results on this is by Kurtz~\cite{Kurtz1978}, who studied continuous-time Markov chains in the context of chemical reaction networks and showed the convergence of the concentrations to the solution of an ordinary differential equation (ODE) as the system size tends to infinity. These results may directly be translated to the Markovian dynamics of interacting agents on an all-to-all coupled network, i.e., on a complete graph~\cite{Niemann2021}. 
Several generalizations of this type of result have emerged, both by employing different dynamics on the network and by considering different interaction networks affecting the agents' dynamics. For an all-to-all coupled network, by allowing the states of the agents to evolve as stochastic differential equations (SDEs), one arrives at what for interacting particle systems is called a McKean--Vlasov equation~\cite{kolokoltsov2010nonlinear}. The mean-field limit of non-Markovian dynamics has been considered e.g.\ in~\cite{duong2018mean,GaRa22}.
As for other interaction structures, one has considered instead of a complete graph also (infinite) lattices~\cite{presutti1983hydrodynamics}, co-evolving graphs~\cite{gkogkas2022graphop}, and in particular random graphs. In the present work, random interaction graphs in conjunction with continuous-time discrete-state random processes are considered.

A recent branch of literature utilizes the concept of graphons~\cite{lovasz2012large,medvedev2014nonlinear} (or graph limits) to derive a large population limit of the dynamics on certain random graphs. For deterministic dynamics, graphon theory is utilized, e.g., in \cite{lee2018consensus} for the analysis of reaching consensus, in \cite{Ayi2020} for graphs with changing edge weights, and in \cite{gkogkas2022graphop} for Kuramoto-type models with adaptive network dynamics.
As for stochastic dynamics, there exist, e.g., works on the stochastic Kuramoto model on Erd{\H o}s--R\'enyi and regular random graphs \cite{delattre2016note} and on particle systems with randomly changing interaction graphs \cite{BhBuWu19}.
In \cite{keliger2022local} the probability distribution of the discrete state of each single node is examined, which in the large population limit yields a description of the process in terms of a partial \review{integro-}differential equation. 
This approach is quite versatile, as the derived limit equation is able to represent a wide range of dynamics. However, this versatility comes with the price of rather high complexity, mathematical intricacies, and also not all random graph models can be modeled by a graphon.
We also note that while \cite{keliger2022local} considers a homogeneous population, where each agent evolves according to analogous (stochastic) rules, our setup will allow for considering a heterogeneous population.

A mean-field limit in the form of a PDE can also be obtained by considering the continuity equation (or transport equation), which is discussed in \cite{Ayi2020, Duteil2022} for deterministic consensus dynamics with time-varying edge weights.
Instead of considering the probability distribution of each node's state, other works describe the large population limit via distributions over other observables, e.g., the number of edges between nodes of certain states~\cite{Decreusefond2012}.

Furthermore, for discrete-state dynamics there are works that exploit symmetries in the network to show that some system states of the global Markov process can be lumped together, in order to obtain a process with fewer states~\cite{banisch2012agent, Simon2010}. While this approach could potentially yield a lower dimensional representation of the process, the finding of (approximately) lumpable states still poses a considerable problem~\cite{Bittracher2021, bittracher2023optimal}.

Another common approach to obtain mean-field type approximations is via so-called ``moment-closure'' methods \cite{miller2014, Porter2016}. Here, equations for the evolution of the frequency of \textit{network motifs} are derived, e.g., the frequency of single nodes in a certain state, the frequency of linked pairs (neighbors) in certain states, or the frequency of triangles in certain states.
These equations are often hierarchical, i.e., the equation for single nodes contains the frequency for pairs, the equation for pairs contains the frequency of triplets, and so on. Hence, a closure of the equations has to be performed to eliminate dependency on higher-order motifs. (Commonly, the equations for pairs are closed, which yields the well-known pair-approximation \cite{Pugliese2009, Peralta2018, Vieira2020}.) This closure introduces an error that generally does not vanish in the large population limit, and hence, these approaches do not guarantee an exact large population limit.

In this work, we prove large population limits in the form of an ordinary differential equation, which we call the \textit{mean-field limit}, for stochastic processes on random graphs.
We consider the case that each node has one of finitely many discrete states and the evolution of a node's state is given by a continuous-time Markov chain.
In this setting, it is of particular interest to observe the shares (or concentrations) of each of the discrete states in the whole system, for example the percentage of infected agents in an epidemiological model.
Hence, we consider these shares as a projection onto the before-mentioned low-dimensional space, and will refer to them as \textit{collective variables}.
However, we also consider collective variables given by finer shares that measure the concentrations of the different states only for a subset of nodes, for instance, the shares in a certain cluster (or community) of nodes.
Our main results are:
\begin{enumerate}
    \item Let a sequence of (random) graphs of increasing size and a dynamical model as described above be given.
    We provide conditions for the choice of collective variables that guarantee that their evolution for the original process converges in probability on finite time intervals to a mean-field ODE (MFE) in the limit of infinitely large networks (Theorem~\ref{thm:main} \review{and Corollary~\ref{cor:quenched}}).
    \item We verify these conditions for multiple random network topologies in section~\ref{sec:voter_model}, in particular for a voter model on Erd\H{o}s--Rényi random graphs, on the stochastic block model, and on random regular graphs. We find that the MFE is a valid approximation for graphs of intermediate density, i.e., where the average degree grows mildly with the graph size---depending on the random graph model, logarithmically or even slower.
    Moreover, we discuss a voter model with heterogeneous population, that is, the dynamical laws of agents can differ across the network.
\end{enumerate}
\review{
While common approaches using graph limits \cite{keliger2022local,Ayi2020} are restricted to random graph models that can be represented by graphons, our high-level convergence theorem (given in section~\ref{sec:main_theorem}) can in principle be applied to arbitrary sequences of random graphs.
For example, random regular graphs exhibit stochastically dependent edges and can thus not be represented via graphons, but they can be handled using the theory presented in this work (cf. section~\ref{subsec:regular} and in particular Corollary~\ref{cor:regular}).
Furthermore, we note that the bounds for the required graph density, which we derive in section~\ref{sec:voter_model}, are consistent with comparable findings from (diluted) graphon approaches~\cite{Bayraktar2020}.
}

The remaining part of this paper is structured as follows.
In section \ref{sec:setup} we give a precise definition of the model and introduce the low dimensional projection (collective variable) we consider.
In section \ref{sec:main_theorem} we prove that the large population limit of the projected dynamics is given by a mean-field ODE, if certain conditions are fulfilled.
We verify these conditions for several examples in section \ref{sec:voter_model} and conclude this paper in section \ref{sec:conclusion}.

\section{Prerequisites} \label{sec:setup}

\paragraph{Notation}
For an integer $N \in \mathbb{N}$, we define $[N]:=\{1,\dots,N\}$.
Furthermore, bold symbols always refer to random variables.
The symbol $\smash{\toprob}$ denotes convergence of random variables in probability.
The symbol $\lVert \cdot \rVert$ denotes any appropriate vector norm.
For two functions $f,g :\mathbb{R} \to \mathbb{R}_{> 0}$ we denote the asymptotic dominance of $f$ over $g$ by
\begin{equation}
    f(x) \gg g(x) \quad :\Leftrightarrow \quad f(x) = \omega(g(x)) \quad :\Leftrightarrow \quad \lim_{x \to \infty} \frac{f(x)}{g(x)} = \infty 
\end{equation}
and asymptotic dominance of $g$ over $f$ by
\begin{equation}
    f(x) \ll g(x) \quad :\Leftrightarrow \quad f(x) = o(g(x)) \quad  :\Leftrightarrow \quad \lim_{x \to \infty} \frac{f(x)}{g(x)} = 0.
\end{equation}

\paragraph{The model}
We consider a simple undirected graph $G = (V, E)$ of size $|V|=N$. 
Without loss of generality, we set $V = [N]$.
Additionally, each node $i\in [N]$ has one of $M \in \mathbb{N}$ discrete states, i.e., $x_i \in [M]$.
We denote the complete system state as $x = (x_1, \dots, x_N) \in [M]^N$.

Each node $i$ changes its state over time due to a continuous-time Markov chain with transition rate matrix
$Q_i^G(x)$, $Q_i^G:[M]^N \to \R^{M\times M}.$
For $m \neq n$, the $(m,n)$-th entry of the rate matrix, $(Q_i^G(x))_{m,n} \geq 0$, specifies at which rate node $i$ transitions from state $m$ to state~$n$. The diagonal entries are such that each row sums to~0.
Note that the transition rates $Q_i^G(x)$ may depend on the full system state $x$ and on the graph~$G$. Moreover, each node $i$ may be subject to a different function $Q_i^G$ determining the transition rates.
We denote the stochastic process referring to the state of node $i$ at time $t$ as $\rv{x}_i(t)$, and the full process as $\rv{x}(t) = \big(\rv{x}_1(t),\dots,\rv{x}_N(t)\big).$
The generator $Q^G$ of the full process $\rv{x}(t)$ is given by
\begin{align}
    Q^G(x,y) = \begin{cases}
    \big(Q_i^G(x)\big)_{x_i, y_i} & \exists i \in [N]\ \forall j \neq i: x_j = y_j\\
    0 & \text{else}
    \end{cases}
\end{align}
and specifies the rate of transitioning to a state $y \in [M]^N$ when starting in a different state $x \in [M]^N$.
Note that this process has the following interpretation: Given $\rv{x}(t) = x$, the processes $\rv{x}_i(s)$, $i=1,\dots,N$, are independent Markov jump processes with rate matrix $Q^G_i(x)$ for $s>t$ as long as no jump takes place, i.e., $\rv{x}(r) = x$ for all~$t<r<s$. When a jump occurs (almost surely no two jumps occur simultaneously), the state $x$ is re-set, with it the rate matrices $Q^G_i(x)$, and the independent processes in the nodes start over with initial conditions dictated by~$x$.
This interpretation gives rise to the common Gillespie algorithm \cite{Gillespie1977} for generating individual realizations of the process: draw the duration until the next event from the exponential distribution corresponding to the accumulated rate of events, then draw the actual event with probabilities proportional to the individual events' rates.

\begin{example} \label{example:common_models}
Common examples of the model introduced above are the following:
\begin{enumerate}
    \item \textbf{SI-model}: In epidemiology, the SI-model \cite{Kiss2017} describes the spreading of a disease on a network. The state of a node is either susceptible (S) or infectious (I). Let $d_{i,I}^G(x) \in \mathbb{N}_0$ denote the number of infectious nodes in the neighborhood of node $i$ when in system state~$x$. Then the rate of the susceptible node $i$ becoming infectious is proportional to $d_{i,I}^G(x)$, and infectious nodes become susceptible again at a constant rate, i.e.,
    \begin{align}
    \label{eq:SIratematrix}
        Q_i^G(x) = \begin{pmatrix}
        - \alpha\, d_{i,I}^G(x) & \alpha\, d_{i,I}^G(x) \\
        \beta & - \beta
        \end{pmatrix}
    \end{align}
    with parameters $\alpha, \beta > 0$.
    
    \item \textbf{Majority rule model}: In the majority rule model of opinion dynamics \cite{Nguyen2020, Porter2016} a node $i$ with state $m$ can only adopt a different state $n$ if a majority of its neighbors have that state $n$. Let the indicator variable $\delta_{i,n}^G(x) \in \{0,1\}$ denote whether more than 50\% of agent $i$'s neighbors have state $n$. Then
    \begin{align}
        \Big(Q_i^G(x)\Big)_{m,n} = \alpha\ \delta_{i,n}^G(x)
    \end{align}
    with parameter $\alpha > 0$. We remark that a deterministic discrete time majority rule model on random graphs is considered in~\cite{sah2021majority}, where a fast convergence to consensus is shown.
    
\item \textbf{Voter model}: In the so-called ``voter model'', originally introduced in \cite{holley1975ergodic}, the state of a node refers to its opinion about some issue and changes stochastically based on the neighborhood of the node. There are many variations of the voter model with slightly different dynamical rules or other modifications, see e.g.~\cite{Moretti2013, Carro2016, Moreira2015, Khalil2018, Peralta2018, Vieira2020}.
In the model that we consider later, a node's rate to change to some state $n$ other than its current state scales linearly with the fraction of nodes in its neighborhood that have state~$n$. We discuss the model in detail in section \ref{sec:voter_model}.
\end{enumerate}
\end{example}

\paragraph{Random graphs}
In the following, we also examine Markov jump processes on random graphs.
We define a random graph $\rv{G}$ on $N$ nodes as a random variable with values in the set of all $2^{N(N-1)/2}$ many simple graphs on nodes $V = [N]$.
In this setting, the stochastic process $\rv{x}(t)$ depends on both the selection of a random graph according to $\rv{G}$ and the stochastic transitions of the Markov jump process.
More precisely, a realization of $\rv{x}(t)$ is given by first sampling a graph $G$ from $\rv{G}$, initializing the node states, and then letting the Markov jump process run on $G$.

\paragraph{Classes and collective variables}
The main result of this work provides conditions under which a low-dimensional representation of the dynamics described above converges to a mean-field ODE in the large population limit.
The map which projects the system state $x$ onto this reduced space is called a collective variable.
We consider a special type of collective variables that measure the concentration of each discrete node state within certain subsets of the nodes. We define these subsets as follows.
We assign each node to one of $K \in \mathbb{N}$ classes and denote the classification of node $i$ as $s_i \in [K]$. Note that the class $s_i$ of node $i$ is fixed and does not depend on the realization of the random graph $\rv{G}$ or on time.
We call the tuple $(x_i, s_i)$ the \textit{extended state} of node~$i$.
Finally, the collective variable $C:[M]^N \to \R^{MK}$ is defined by measuring the shares of each extended state in each of these subsets, i.e.,
\begin{align} \label{eq:C}
    C(x) = \Big(C_{(m, k)}(x)\Big)_{m \in [M], k \in [K]} , \qquad C_{(m, k)}(x) := \frac{\# \{i \in [N]: (x_i, s_i) = (m, k) \}}{N}.
\end{align}
It is well-known that mean-field approximations, i.e., expressing the dynamics only in terms of concentrations, work best if nodes (or particles in physical literature) are at least somewhat indistinguishable and interchangeable \cite{Porter2016,Ayi2020}. Hence, it is reasonable in our setting to group nodes together in classes if they have similar traits and thus may become indistinguishable in the large population limit. We provide some examples below.

\begin{figure}
\centering
\includegraphics[width=.8\textwidth]{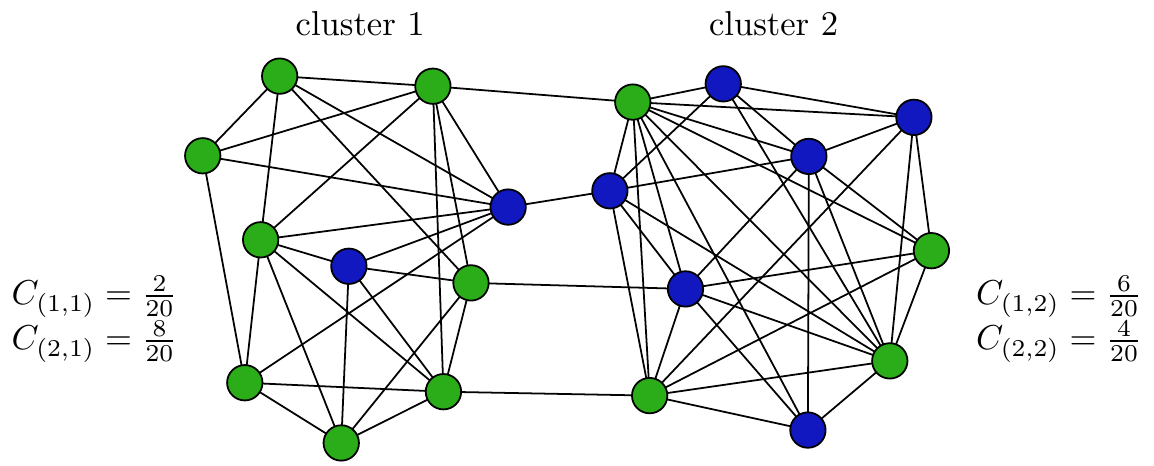}
\caption{Example graph $G$ of size $N=20$, sampled from a stochastic block model with two clusters (see section \ref{subsec:sbm} for more details). Each node has one of two states, state $1$ is indicated by blue color and state $2$ by green color. We assign a node to class $k$ if it is located in cluster $k$, $k=1,2$. The collective variable $C(x)$ measures the shares of the two states in each cluster.}
\label{fig:illustration}
\end{figure}

\begin{example} \label{example:common_species}
Common examples for the choice of classes:
\begin{enumerate}
    \item In the case $K=1$, the extended states are the states itself, $(x_i, s_i) = (x_i, 1) \cong x_i$. The collective variable $C(x)$ measures the share of each state in the system. These collective variables are commonly discussed in mean-field literature, e.g.~\cite{slanina2003analytical}.
    A necessary and sufficient condition on the random graph sequence for the mean field limit to hold has been derived in \cite{keliger2021_quasi} for processes where the transition rates $Q_i$ are affine-linear in the so-called neighborhood vector (an example would be given by equation~\eqref{eq:SIratematrix}). We note that the continuous-time noisy voter model that we consider later (cf.\ section~\ref{subsec:CNVM}) does not fall in this category of processes and neither do so-called ``complex contagion'' models in which the infection rates are nonlinear functions of the infection prevalence in the node's neighborhood.
    \item If the random graph exhibits a fixed modular structure with $K$ communities or clusters, it is a natural choice to measure the shares of the states in each cluster separately, i.e., a node located in cluster $k$ is assigned to class $k$. Hence, the extended state $(x_i, s_i) = (m, k)$ refers to a node with state $x_i = m$ located in cluster $k$.
    We provide an example in Figure~\ref{fig:illustration}.
    Differentiating nodes by their community is frequently used in the literature, see e.g.~\cite{Khalil2018, niemann2021a, Helfmann2021a}.
    The well-known stochastic block model, which we discuss in more detail in section \ref{subsec:sbm}, generates random graphs that exhibit this clustered structure. 
    \item If nodes differ by their transition rate matrices $Q_i^G$, i.e., the population is heterogeneous, it is reasonable to assign them to different classes. As an example, there could be very active nodes (large transition rates) and rather inactive nodes (small or zero transition rates, often referred to as ``zealots'' in literature \cite{huo2019zealot, Khalil2018, Chinellato2015}). We discuss the case of a heterogeneous population in section \ref{subsec:hetero}.
    \item We may want to differentiate between nodes by their position on the graph. For instance, we could construct classes based on node degrees or centrality measures. An ``influencer'' class could consist of nodes with high node degrees, while the ``follower'' class has lower degrees. Note that in this case the class $s_i$ of node $i$ depends on the realization of the random graph $\rv{G}$, contrary to our previous assumption. We can remedy this issue by defining a function $s$ such that $s_i = s(i, G)$ and adapting the collective variable accordingly, i.e.,
    \begin{align}
        C^G_{(m, k)}(x) := \frac{\# \{i \in [N]: (x_i, s(i, G)) = (m, k) \}}{N}.
    \end{align}
    The main theorem for convergence in the large population limit, which we will discuss in section~\ref{sec:main_theorem}, also applies for this slight extension of the class framework. However, as the examples that we discuss later in section \ref{sec:voter_model} all employ fixed classes $s_i$, we will not consider this extension further in this work.
\end{enumerate}
\end{example}

Due to the dynamics, a node with extended state $(m,k)$ may transition to any other state $n \neq m$, such that after the transition it has extended state $(n,k)$.
We refer to this transition as $(m,k)\to n$ and there are $MK (M-1)$ transitions in total.
Each transition has an associated state-change vector $v_{(m,k)\to n} \in \mathbb{Z}^{MK}$ and a propensity function $\alpha^G_{(m,k)\to n}:[M]^N \to \R_{\geq 0}$.
The state-change vector
\begin{align}
    v_{(m,k)\to n} := e_{(n,k)} - e_{(m,k)},
\end{align}
where $e_{(m,k)}$ denotes the $(m,k)$-th unit-vector, describes the changes in extended state populations due to the associated transition $(m,k)\to n$, i.e., there is one less node in extended state $(m,k)$ and one more in~$(n, k)$.
The propensity function $\alpha^G_{(m,k)\to n}$ measures the cumulative rate of the transition $(m,k)\to n$, i.e., the sum of the transition rates of all agents with extended state $(m,k)$ to state~$n$
\begin{align}
    \alpha^G_{(m,k)\to n}(x) := \sum_{i \in [N]: (x_i,s_i) = (m,k)} \big(Q_i^{G}(x)\big)_{m, n}. \label{eq:propensity_functions}
\end{align}

In the following, we abbreviate the summation over all $MK (M-1)$ transitions with the symbol
\begin{align}
    \sum_{(m,k)\to n} := \sum_{m=1}^M \sum_{k=1}^K \sum_{\substack{n=1 \\ n\neq m}}^M,
\end{align}
and analogously for the maximum over all transitions: $\max_{(m,k)\to n}$.

\section{Conditions for convergence in the large population limit} \label{sec:main_theorem}

Let the extended states $(m, k)$ and collective variables $C(x)$ be as defined in the previous section.
We assume that the classes are chosen such that the collective variables capture the most important dynamical information. This statement is made more precise later in the conditions of Theorem~\ref{thm:main}. In a broad sense we demand that the propensities \eqref{eq:propensity_functions} of the transitions can be approximated well by using only the reduced information $C(x)$ of the state, i.e.,
there exist propensity functions $\tilde{\alpha}_{(m,k)\to n}: \mathbb{R}^{MK} \to \R$ with
\begin{align}
    \frac{1}{N} \alpha_{(m,k)\to n}^{G}(x)  \approx \tilde{\alpha}_{(m,k)\to n}\big( C(x) \big) \quad \forall x\in [M]^N.
    \label{eq:propensity}
\end{align}
We assume that all $\tilde{\alpha}_{(m,k)\to n}$ are Lipschitz continuous.
Existence of such an appropriate choice of classes and reduced propensity functions $\tilde{\alpha}_{(m,k)\to n}$ for a given dynamical system on a certain network is not clear, and finding them is no trivial task.
However, if classes and propensity functions can be found such that the approximation \eqref{eq:propensity} becomes exact in the large population limit, then there exists a mean-field ODE describing the projected system state $C(x)$, which we show in Theorem \ref{thm:main}.

In order to specify the large population limit, we consider a sequence of random graphs $(\rv{G}_\ell)_{\ell \in \mathbb{N}}$, such that $\rv{G}_\ell$ has $N_\ell$ nodes and the sequence $(N_\ell)_{\ell \in \mathbb{N}}$ is strictly increasing.
Furthermore, let $s_i^\ell \in [K]$ denote the class of node $i$ of the random graph $\rv{G}_\ell$ and define the collective variables
\begin{align}
    C^\ell_{(m, k)}(x) := \frac{\# \{i \in [N_\ell]: (x_i, s^\ell_i) = (m, k) \}}{N_\ell}.
\end{align}
Let $\rv{x}^\ell(t)$ denote the stochastic jump process on the random graph $\rv{G}_\ell$, and
\begin{align}
    \rv{c}^\ell(t) := C^{\ell}(\rv{x}^\ell(t))
    \label{eq:definition_c}
\end{align}
the projected process.
In order to quantify the approximation \eqref{eq:propensity}, we define the difference
\begin{align}\label{Delta}
    \Delta^{\rv{G}_\ell}(x) := \max_{(m,k)\to n} \Big\lvert \frac{1}{N_\ell} \alpha_{(m,k)\to n}^{\rv{G}_\ell}(x) - \tilde{\alpha}_{(m,k)\to n}\big( C^\ell(x) \big) \Big\rvert
\end{align}
for $x \in [M]^{N_\ell}$.
Moreover, we assume that the transition rate matrices $Q_i^{G}(x)$ have bounded entries, i.e., there is a bound $B > 0$, such that for any graph $G$ on any number of nodes we have
\begin{align}
    Q_i^{G}(x) < B \quad \text{elementwise for any system state $x$ and all $i\in [N]$.}
    \label{eq:Q_bound_B}
\end{align}

The mean-field ODE, which we will refer to as the mean-field equation (MFE), is given by
\begin{align}
    \frac{\diff}{\diff t} c(t) = \sum_{(m,k)\to n} \tilde{\alpha}_{(m,k)\to n}(c(t))\ v_{(m,k)\to n} =: F\big(c(t)\big),
    \label{eq:MFE}
\end{align}
where $c:\R \to \R^{MK}$.
The infinitesimal change in $c$ is characterized by the propensities of each transition $\tilde{\alpha}_{(m,k)\to n}$ times their effect on the extended state populations~$v_{(m,k)\to n}$.
Due to the assumption that the $\tilde{\alpha}_{(m,k)\to n}$ are Lipschitz continuous, the mean-field ODE has a unique solution, given initial condition~$c(0) = c_0$.

The following theorem is inspired by work about the concentration of Markov processes by Kurtz \cite{Kurtz1978}, which can be applied to a sequence of complete graphs increasing in size. The proof of the following theorem generalizes the proof of Kurtz's result presented in~\cite[Thm. A.14]{Winkelmann2020}, and is relying on combining the law of large numbers and Gronwall's lemma.
\begin{theorem} \label{thm:main}
Assume that for all $\varepsilon > 0$ there exists a function $f_\varepsilon:\mathbb{N}\to \R_{\geq 0}$ with $\lim_{\ell \to \infty} f_\varepsilon(\ell) = 0$ such that
\begin{align}
    \forall \ell \in \mathbb{N}:\ \mathbb{P}\Big(\max_{x\in [M]^{N_\ell}} \Delta^{\rv{G}_\ell}(x) \geq \varepsilon \Big) \leq f_\varepsilon(\ell).
\end{align}
Furthermore, let there be initial conditions $\rv{x}^\ell(0)$, such that $\rv{c}^\ell(0) \toprob c_0 \in \R^{MK}$ as $\ell \to \infty$.
Let $c(\cdot)$ denote the solution of the mean-field equation \eqref{eq:MFE} with initial condition $c(0) = c_0$.
Then 
\begin{align}
    \forall t \geq 0:\ \sup_{0 \leq s \leq t} \big\lVert \rv{c}^\ell(s) - c(s) \big\rVert \underset{\ell \to \infty}{\toprob} 0.
\end{align}
\end{theorem}

\begin{proof}
See \cref{sec:appendix_proof_main}.
\end{proof}

\review{
\begin{remark}[Rate of convergence]
    In the proof of the previous theorem, we derive the following bound (cf. equation~\eqref{eq:approx_quality}) which indicates the main factors controlling the rate of convergence as $\ell \to \infty$:
    \begin{align}
    \lVert \rv{c}^\ell(t) - c(t) \rVert \leq \Big(\lVert \rv{c}^\ell(0) - c(0) \rVert + \hat{\rv{\delta}}_\ell(t) + \rv{\tilde{\delta}}_\ell(t) \Big) \exp(L t).
\end{align}
In this remark we roughly outline how the above bound behaves. For a detailed definition of the occurring symbols see \cref{sec:appendix_proof_main}.
First of all, we note that the factor $\exp(L t)$ implies that the deviation between the stochastic process $\rv{c}^\ell(t)$ and the mean-field solution $c(t)$ may increase over time for a fixed $\ell$. The rate $L$ of this deterioration is proportional to the Lipschitz constants of the propensities $\tilde{\alpha}_{(m,k)\to n}$.
Hence, for practical purposes, if a good match between model and mean-field solution is required for a longer time, the network size $N_\ell$ may have to be increased substantially.
For the three terms inside the brackets, we note the following:
\begin{enumerate}
    \item The difference between initial conditions $\lVert \rv{c}^\ell(0) - c(0) \rVert$ is typically not a limiting factor. For example, for any target initial condition $c_0 \in [0,1]^{MK}$, we can find (deterministic) initial conditions $\rv{c}^\ell(0)$ on the graph with an error $\lVert \rv{c}^\ell(0) - c_0 \rVert = O(N_\ell^{-1})$.
    \item The convergence of the second term $\hat{\rv{\delta}}_\ell(t) \to 0$ is essentially given by the law of large numbers, applied to normalized and centered Poisson processes.
    Thus, the rate of convergence of $\sqrt{N_\ell}^{-1}$ is dictated by the central limit theorem.
    For a detailed derivation, see Remark \ref{remark_appendix_convergence}.
    \item Finally, the term $\rv{\tilde{\delta}}_\ell(t) \geq 0$ captures the influence of the random graph. From \eqref{eq:tilde_delta_def} and Lemma \ref{lemma:conv_prob}, it follows $\mathbb{E}[\rv{\tilde{\delta}}_\ell(t)] = O(f_\varepsilon(\ell))$ for all~$\varepsilon$. Hence, the rate at which this term goes to 0 is determined by how fast $\Delta^{\rv{G}_\ell}$ becomes small.
\end{enumerate}
We can conclude that the rate of convergence is the slower one of $\sqrt{N_\ell}^{-1}$ (due to point 2.) and $f_\varepsilon(\ell)$ (due to point 3.).
\end{remark}
}

\begin{remark} \label{remark:main_deterministic}
The case of a sequence of deterministic graphs $(G_\ell)_{\ell \in \mathbb{N}}$ (as opposed to random graphs) is also contained in the previous theorem. In this case, the condition of the theorem collapses to
\begin{align}
    \max_{x \in [M]^{N_\ell}} \Delta^{G_\ell}(x) \overset{\ell \to \infty}{\longrightarrow} 0.
\end{align}
\end{remark}

\review{
\begin{corollary}[Quenched result] \label{cor:quenched}
    If the function $f_\varepsilon\ $ from Theorem \ref{thm:main} additionally satisfies
    \begin{equation}
        \forall \varepsilon > 0:\ \sum_{\ell = 1}^\infty f_\varepsilon(\ell) < \infty,
    \end{equation}
    the convergence to the mean-field limit also holds for almost all realizations of the sequence of random graphs.
    More precisely, if $\rv{c}^\ell_{G}(s)$ denotes the stochastic process given by the network dynamics on the fixed graph $G$, then for almost all realizations $(G_\ell)_\ell$ of the sequence of random graphs $(\rv{G}_\ell)_\ell$ we have 
    \begin{equation}
    \forall t \geq 0:\ \sup_{0 \leq s \leq t} \big\lVert \rv{c}^\ell_{G_\ell}(s) - c(s) \big\rVert \underset{\ell \to \infty}{\toprob} 0. \label{eq:main_convergence_quenched}
    \end{equation}
    This stronger result is often called the quenched result, whereas the previous Theorem \ref{thm:main}, in which $\rv{c}^\ell(s)$ contains both the random graph and the random dynamics, describes the annealed result.
\end{corollary}
\begin{proof}
    From the Borel--Cantelli lemma it follows that with probability $1$ the event
    \begin{align}
        \left\{ \max_{x\in [M]^{N_\ell}} \Delta^{\rv{G}_\ell}(x) \geq \varepsilon \right\}
    \end{align}
    only occurs for finitely many $\ell$.
    Due to the arbitrary choice of $\varepsilon > 0$, we have
    \begin{align}
        \mathbb{P}\left( \max_{x\in [M]^{N_\ell}} \Delta^{\rv{G}_\ell}(x) \underset{\ell \to \infty}{\longrightarrow} 0 \right) = 1.
    \end{align}
    Hence, we can apply Remark \ref{remark:main_deterministic} for almost all realizations of the sequence of random graphs.
\end{proof}
}

\section{Large population limits for a voter model} \label{sec:voter_model}

In this section we analyze the large population limit of the continuous-time noisy voter model (CNVM) on Erd\H{o}s--Rényi random graphs (cf.\ section~\ref{subsec:ER}), the stochastic block model (cf.\ section~\ref{subsec:sbm}), and uniformly random regular graphs (cf.\ section~\ref{subsec:regular}) by verifying that the conditions of Theorem \ref{thm:main} hold. Moreover, we also consider a heterogeneous population, where the agents may have different transition rates, in section~\ref{subsec:hetero}.
We start with introducing the CNVM of our interest.

\subsection{The continuous-time noisy voter model (CNVM)}
\label{subsec:CNVM}

The CNVM originates from opinion dynamics, but is connected to many other applications like epidemiology \cite{Kiss2017,Porter2016}, genetics \cite{moran1958}, and statistical mechanics \cite{BINDER1972}. Hence, we refer to the state $x_i$ of node $i$ as its opinion.
Moreover, we use the terms ``node'' and ``agent'' synonymously.
The CNVM can be described by the following procedure:
\begin{enumerate}
    \item Each agent is regularly influenced by its neighbors. For all agents (that have at least one neighbor), these influencing events happen randomly with rate $r_0>0$ and independently from others.
    \item If an event is triggered for agent $i$, a random neighbor is chosen. Hence, the probability that agent $i$ is influenced by a neighbor of opinion $n$ is given by $\smash{ d_{i,n}^G(x) / d_{i}^G }$, where $\smash{ d_{i,n}^G(x) }$ denotes the number of neighbors of agent $i$ with opinion $n$, and $\smash{ d_{i}^G }$ the node degree of agent~$i$.
    \item After being influenced, agent $i$ decides whether or not to adopt the opinion of the contacted neighbor. Let $p_{m,n}\in [0,1]$ denote the probability that an agent of opinion $m$, after being influenced by a neighbor of opinion $n$, adopts the opinion~$n$.
    \item Additionally, agent $i$ transitions from its current opinion $m$ to a different opinion $n$ at a constant rate $\tilde{r}_{m,n} \geq 0$ that does not depend on its neighbors. (Hence, the $\tilde{r}_{m,n}$ can be thought of as noise intensity.)
\end{enumerate}
Thus, the entries of the transition rate matrices $Q_i^G(x) \in \mathbb{R}^{M\times M}$ are given by
\begin{align}
    \big(Q_i^G(x)\big)_{m,n} &= r_0 \frac{d_{i,n}^G(x)}{d_{i}^G} p_{m,n} + \tilde{r}_{m,n} \\
    &= r_{m, n} \frac{d_{i,n}^G(x)}{d_{i}^G} + \tilde{r}_{m,n}, \label{eq:CNVM_rates}
\end{align}
where $r_{m, n} := r_0 p_{m,n}$ and $m\neq n$.
Due to the various parameters, i.e., the number of opinions $M \in \mathbb{N}$, the rates $r_{m,n} \geq 0$, and the rates $\tilde{r}_{m,n} \geq 0$, the CNVM is a fairly general model for a simple spreading process on a graph.

The following example from epidemiology is in its abstract form also a CNVM:
\begin{example}[SIRS model \cite{Volz2007}] \label{example:SIR}
We construct an SIRS (susceptible, infectious, recovered, susceptible) model from the CNVM as follows. If an agent is susceptible, it has a rate of $r_{S,I} > 0$ to become infectious if all its neighbors are infectious. This rate scales linearly with the share of infectious neighbors, e.g., if half of its neighbors are infectious the agent becomes infectious at a rate of $\frac{1}{2} r_{S,I}$.
If an agent is infectious, it takes on average $ 1 / \tilde{r}_{I,R}$ time units until the infection is over and the agent becomes recovered. (To be precise, the event is exponentially distributed with rate $\tilde{r}_{I,R}$.)
Being recovered, an agent becomes susceptible again after $1 / \tilde{r}_{R,S}$ time units on average.
Hence, the rates $r_{m,n}$ and $\tilde{r}_{m,n}$, written in matrix form in the order $(S,I,R)$, are given by
\begin{align}
r = \begin{pmatrix}
- & r_{S,I} & 0\\
0 & - & 0\\
0 & 0 & -
\end{pmatrix}, \quad \tilde{r} = \begin{pmatrix}
- & 0 & 0\\
0 & - & \tilde{r}_{I,R}\\
\tilde{r}_{R,S} & 0 & -
\end{pmatrix}.
\end{align}
\end{example}

In the following sections we derive and prove mean-field limits of the CNVM on different random graphs.

\subsection{Erd\H{o}s--Rényi random graphs} \label{subsec:ER}

The Erd\H{o}s--Rényi (ER) random graph $\rv{G}_{N,p}$, also called binomial random graph, is defined as the random graph where each possible edge appears independently with probability $p>0$. To be precise, we have
\begin{align}
    \mathbb{P}(\rv{G}_{N, p} = G) = p^{\lvert E(G) \rvert} (1 - p)^{\binom{N}{2} - \lvert E(G) \rvert},
\end{align}
where $G$ is a simple graph on vertices $[N]$ and $|E(G)|$ denotes the number of edges in $G$~\cite{Frieze2015}.
\review{We implicitly allow the edge probability $p = p(N)$ to depend on the number of nodes $N$.
It is especially interesting to investigate if a mean-field limit exists depending on the asymptotic behavior of $p$, e.g., depending on how fast $p$ converges to $0$ as $N \to \infty$.}

In this section we show that the CNVM on ER random graphs converges to a mean-field limit with respect to the shares of each opinion, provided the expected node degree $p(N-1)$ grows fast enough with $N$.
Thus, we have $K=1$ and the extended states $(x_i, s_i) = (x_i, 1)$ collapse to just the states~$x_i$.
For easier notation, we refer to the transition $(m, 1) \to n$ as $m \to n$.
Moreover, we consider the sequence of random graphs $\big(\rv{G}_{N,p}\big)_{N\in\mathbb{N}}$, such that we can use the index $N$ instead of $\ell$ from the previous section.

\paragraph{Heuristic derivation of the mean-field equation}
Let us first derive the mean-field ODE in a heuristic manner. Consider for a given graph $G$ the propensity functions
\begin{align}
    \alpha^{G}_{m\to n}(x) &= \sum_{i: x_i = m} \big(Q_i^{{G}}(x)\big)_{m, n} \\
    &= \sum_{i: x_i = m} \Big( r_{m, n}  \frac{d_{i,n}^{G}(x)}{d_{i}^{G}} + \tilde{r}_{m,n} \Big)\label{eq:cnvm_propensity}
\end{align}
\review{describing the cumulative rate at which an opinion change $m\to n$ occurs in the entire graph.}
Because of the homogeneous nature of an ER random graph $\rv{G}$, we expect the share of agents of opinion $n$ in the neighborhood of agent $i$, $d_{i,n}^{\rv{G}}(x) / d_{i}^{\rv{G}}$, to be approximately equal to the share of opinion $n$ in the whole system, $C_n(x) =: c_n$.
Hence, we have
\begin{align}
    \frac{1}{N} \alpha^{\rv{G}}_{m\to n}(x) &\approx \frac{1}{N} \sum_{i: x_i = m} \Big( r_{m, n} c_n + \tilde{r}_{m,n} \Big)  \\
    &= c_m \big(r_{m, n} c_n +  \tilde{r}_{m,n} \big) =: \tilde{\alpha}_{m\to n}(c), \label{eq:cnvm_reduced_propensity}
\end{align}
which yields the following mean-field ODE when inserted into \eqref{eq:MFE}:
\begin{align} \label{eq:MFE_ER}
    \frac{\diff}{\diff t} c(t) = \sum_{m \neq n} c_m(t) \big(r_{m, n} c_n(t) +  \tilde{r}_{m,n} \big) (e_n - e_m),
\end{align}
\review{where the summation is over all pairs $(m,n) \in [M] \times [M]$ with~$m\neq n$.}

Now we show that the propensity functions $\tilde{\alpha}_{m\to n}(c)$ derived above indeed fulfill the conditions of Theorem \ref{thm:main}. We will make use of the following results.

\paragraph{Auxiliary concentration results}

\begin{lemma} \label{lemma:ER_counts}
Let $\rv{G} = \rv{G}_{N,p}$ denote the ER random graph and $x \in [M]^N$ an arbitrary but fixed state. Define the random variable $\rv{E}_{m,n}$ as the number of edges between nodes of opinion $m$ and nodes of opinion $n\neq m$ in $\rv{G}$, according to $x$, i.e., $\rv{E}_{m,n} = \sum_{i:x_i=m} d_{i,n}^{\rv{G}}(x)$.
Then we have the concentration inequality
\begin{align}
    \mathbb{P}(|\rv{E}_{m,n} - c_m c_n N^2 p | \geq \varepsilon) \leq 2 \exp\Big(-\frac{ \varepsilon^2}{3 N^2 p}\Big)
\end{align}
for all $\varepsilon>0$, where~$c := C(x)$.
\end{lemma}
\begin{proof}
    There are $u := c_m c_n N^2$ possible edges between $m$-opinion and $n$-opinion nodes. As every edge in $\rv{G}_{N,p}$ is present with probability $p$ independently of all other edges, it follows that $\rv{E}_{m,n}$ is binomially distributed with $u$ trials and success probability $p$ for each trial.
    In particular, we have $\mathbb{E}[\rv{E}_{m,n}] = u p = c_m c_n N^2 p$, and applying the Chernoff bound (Lemma \ref{lemma:chernoff}) yields
    \begin{align}
    \mathbb{P}(|\rv{E}_{m,n} - c_m c_n N^2 p | \geq \varepsilon) &\leq 2 \exp\Big(- \frac{\varepsilon^2}{3 u p}\Big) \\
    &\leq 2 \exp\Big(-\frac{ \varepsilon^2}{3 N^2 p}\Big),
    \end{align}
    where the last inequality follows from~$u \leq N^2$.
\end{proof}

Furthermore, we can show in a similar fashion (see Lemma~\ref{lemma:A_ER_degrees}) that the node degrees $\smash{ d_i^{\rv{G}_{N,p}} }$ are concentrated around $N p$, i.e., for $\varepsilon > 0$ and $i \in [N]$ we have
\begin{align}
    \mathbb{P}\left(\big|d_i^{\rv{G}_{N,p}} - Np\big| \geq \varepsilon Np \right) \leq 2 \exp\Big( -\frac{\varepsilon^2 N p}{3} + \frac{2 \varepsilon}{3}\Big). \label{eq:ER_degrees}
\end{align}
Now we can verify the conditions of Theorem~\ref{thm:main}:
\begin{proposition} \label{prop:ER_f}
Let $\rv{G}_{N,p}$ denote the ER random graph and $\hat{r} := \max_{m\neq n} r_{m,n}$.
For all $\varepsilon \in (0, \hat{r})$  we have that
\begin{align}
    \forall N \in \mathbb{N}:\ \mathbb{P}\Big(\max_{x\in [M]^{N}} \Delta^{\rv{G}_{N, p}}(x) \geq \varepsilon \Big) \leq f_\varepsilon(N),
\end{align}
where
\begin{align} \label{f_eps}
    f_\varepsilon(N) := 4 M^{N+2} \exp\Big(-\frac{1}{12} N^2 p \Big(\frac{\varepsilon}{\hat{r}} - \frac{\varepsilon^2}{\hat{r}^2}\Big)^2 \Big) + 2 N \exp\Big(-N \frac{\varepsilon^2 p}{12 \hat{r}} + \frac{ \varepsilon}{3\hat{r}}\Big).
\end{align}
\end{proposition}
\begin{proof}[Proof sketch.]
The detailed proof can be found in \cref{sec:appendix_proof_ER_f}. The main idea of the proof is as follows.
We fix any $N \in \mathbb{N}$ and denote $\rv{G} := \rv{G}_{N,p}$.
By inserting the propensity functions \eqref{eq:cnvm_propensity} and \eqref{eq:cnvm_reduced_propensity} into \eqref{Delta}, it follows that
\begin{align}
    \Delta^{\rv{G}}(x) = \max_{m\neq n}\  r_{m,n} \Big\lvert \frac{1}{N} \sum_{i : x_{i} = m} \frac{d^{\rv{G}}_{i,n}(x)}{d^{\rv{G}}_i} - C_m(x) C_n(x) \Big\rvert.
\end{align}
    Let $\delta \in (0,1)$ and define the random events
    \begin{align}
        \mathcal{A} &:= \Big\{ \max_{x\in [M]^{N}} \Delta^{\rv{G}}(x) \geq \varepsilon \Big\} = \Big\{ \max_{x\in [M]^{N}} \max_{m \neq n}\ r_{m,n} \Big\lvert \frac{1}{N} \sum_{i : x_{i} = m} \frac{d^{\rv{G}}_{i,n}(x)}{d^{\rv{G}}_i} - C_m(x) C_n(x) \Big\rvert \geq \varepsilon \Big\},\\
        \mathcal{B} &:= \Big\{ \forall i: (1-\delta)Np \leq d^{\rv{G}}_i \leq (1+\delta)Np\Big\}.
    \end{align}
    We need to bound~$\mathbb{P}(\mathcal{A})$. Note that, due to the concentrated node degrees in ER random graphs, the probability of event $\mathcal{B}$ is large. To be precise, from equation~\eqref{eq:ER_degrees} and the union bound, it follows that
    \begin{align}
        \mathbb{P}(\mathcal{B}^C) \leq 2 N \exp\Big(\frac{1}{3} (-\delta^2 N p + 2 \delta)\Big),
    \end{align}
    where $\mathcal{B}^C$ denotes the complement of $\mathcal{B}$.
    Using Lemma \ref{lemma:ER_counts} and an appropriate choice of $\delta$, we can show that
    \begin{align}
        \mathbb{P}(\mathcal{A} \cap \mathcal{B}) \leq 4 M^{N+2} \exp\Big(-\frac{1}{12} N^2 p \Big(\frac{\varepsilon}{\hat{r}} - \frac{\varepsilon^2}{\hat{r}^2}\Big)^2 \Big).
    \end{align}
    All in all, this yields
    \begin{align}
         \mathbb{P}(\mathcal{A}) &\leq \mathbb{P}(\mathcal{A} \cap \mathcal{B}) + \mathbb{P}(\mathcal{B}^C)\\
        &\leq 4 M^{N+2} \exp\Big(-\frac{1}{12} N^2 p \Big(\frac{\varepsilon}{\hat{r}} - \frac{\varepsilon^2}{\hat{r}^2}\Big)^2 \Big) + 2 N \exp\Big(-N \frac{\varepsilon^2 p}{12 \hat{r}} + \frac{ \varepsilon}{3\hat{r}}\Big),
    \end{align}
    from which the proposition follows.
\end{proof}

\paragraph{Large population limit}

By examining the bounding function $f_\varepsilon$ that we have derived in the previous proposition, we can conclude edge densities $p = p(N)$ for which the mean-field result holds. The following theorem states that ER random graphs of intermediate density are sufficient to obtain the mean-field limit. Interestingly, the derived threshold for the edge density $p$ is exactly the sharp threshold that yields (asymptotically almost surely) connectedness of $\rv{G}_{N,p}$ \cite{Frieze2015}.
\begin{theorem} \label{thm:ER_convergence}
Let the edge probability $p=p(N)$ of the ER random graph $\rv{G}_{N, p}$ be a function of the number of vertices~$N$.
If $p$ dominates $\log(N)/N$ asymptotically, i.e.,
\begin{align}
    p = \omega\Big(\frac{\log N}{N}\Big) \quad \text{as}\ N \to \infty, 
\end{align}
then the dynamics of the opinion shares in the CNVM converges to a mean-field limit as $N \to \infty$, \review{in the sense of both Theorem~\ref{thm:main} (annealed result) and Corollary~\ref{cor:quenched} (quenched result)}. The mean-field solution satisfies the ODE
\begin{align}
    \frac{\diff}{\diff t} c(t) = \sum_{m \neq n} c_m \big(r_{m, n} c_n +  \tilde{r}_{m,n} \big) (e_n - e_m). \label{eq:mfe_ER_2}
\end{align}
\end{theorem}
\begin{proof}
    In Proposition \ref{prop:ER_f} we derived the function 
    \begin{align}
        f_\varepsilon(N) = 4 M^{N+2} \exp\Big(-\frac{1}{12} N^2 p \Big(\frac{\varepsilon}{\hat{r}} - \frac{\varepsilon^2}{\hat{r}^2}\Big)^2 \Big) + 2 N \exp\Big(-N \frac{\varepsilon^2 p}{12 \hat{r}} + \frac{ \varepsilon}{3\hat{r}}\Big) \label{eq:f_epsilon_ER}
    \end{align}
    as a bound for $\mathbb{P}\big(\max_{x} \Delta^{\rv{G}_N}(x) \geq \varepsilon \big)$. In order to apply Theorem \ref{thm:main}, we have to make sure that $f_\varepsilon(N) \to 0$ as $N \to \infty$.
    For the right term in $f_\varepsilon$ to converge to $0$, it is necessary that $Np\varepsilon^2$ dominates $\log N$ for all $\varepsilon > 0$, which is given for $p = \omega\big(\frac{\log N}{N}\big)$.
    Similarly, for the left term in $f_\varepsilon$ to converge to $0$, it is necessary that $N^2 p$ dominates $N$, which is less restrictive and also true for~$\smash{ p = \omega\big(\frac{\log N}{N}\big) }$.

    \review{Moreover, neglecting constants, the right term of $f_\varepsilon$ satisfies $N \exp(-N p) \ll N^{-2}$, from which the condition $\sum_N f_\varepsilon(N) < \infty$ of Corollary~\ref{cor:quenched} follows.}
\end{proof}

\begin{example} \label{example:ER}
We provide numerical results of two example models in Figure \ref{fig:concentration_ER}, which indicate how the derived mean-field solution becomes a good approximation of the stochastic process $\rv{c}(t)$ for large numbers of agents $N$. 
\review{For every numerical sample of the process, we generate a new random graph and then simulate the CNVM on that graph using Gillespie's stochastic simulation algorithm \cite{Gillespie1977}.}
In the first example (cf.\ Figure~\ref{fig:concentration_ER_a}) we examine the CNVM with $M=2$ opinions, initial conditions $c(0) = (0.2, 0.8)$, and rate constants
\begin{align}
    r = \begin{pmatrix}
         - & 0.99\\
         1 & - \\
         \end{pmatrix},
    \quad \tilde{r} = \begin{pmatrix}
         - & 0.01\\
         0.01 & - \\
         \end{pmatrix}.
\end{align}
In the second example (cf.\ Figure~\ref{fig:concentration_ER_b}) there are $M=3$ opinions, the initial condition is $c(0) = (0.2, 0.5, 0.3)$, and rate constants are
\begin{align}
    r = \begin{pmatrix}
         - & 0.8 & 0.2\\
         0.2 & - & 0.8\\
         0.8 & 0.2 & - \\
         \end{pmatrix},
    \quad \tilde{r}_{m,n} = 0.01 \text{ for all } m\neq n.
\end{align}
For both examples the edge density was set to~$p=0.01$. 
Note that if the number of agents is too small, the mean-field equation fails to approximate the dynamics well because of the high variance of $\rv{c}(t)$, and also the mean of $\rv{c}(t)$ may not be close to the mean-field solution.
\review{As we let the number of agents increase, the variance of the process decreases and the mean moves closer to the mean-field solution, see Figure~\ref{fig:concentration_ER_a}.}
Moreover, note that the quality of the approximation of $\rv{c}(t)$ by the mean-field limit may deteriorate over time, as indicated by equation \eqref{eq:approx_quality}, which can be seen in Figure~\ref{fig:concentration_ER_b}.

\begin{figure}
     \centering
     \begin{subfigure}[t]{0.49\textwidth}
         \centering
         \includegraphics[width=\textwidth]{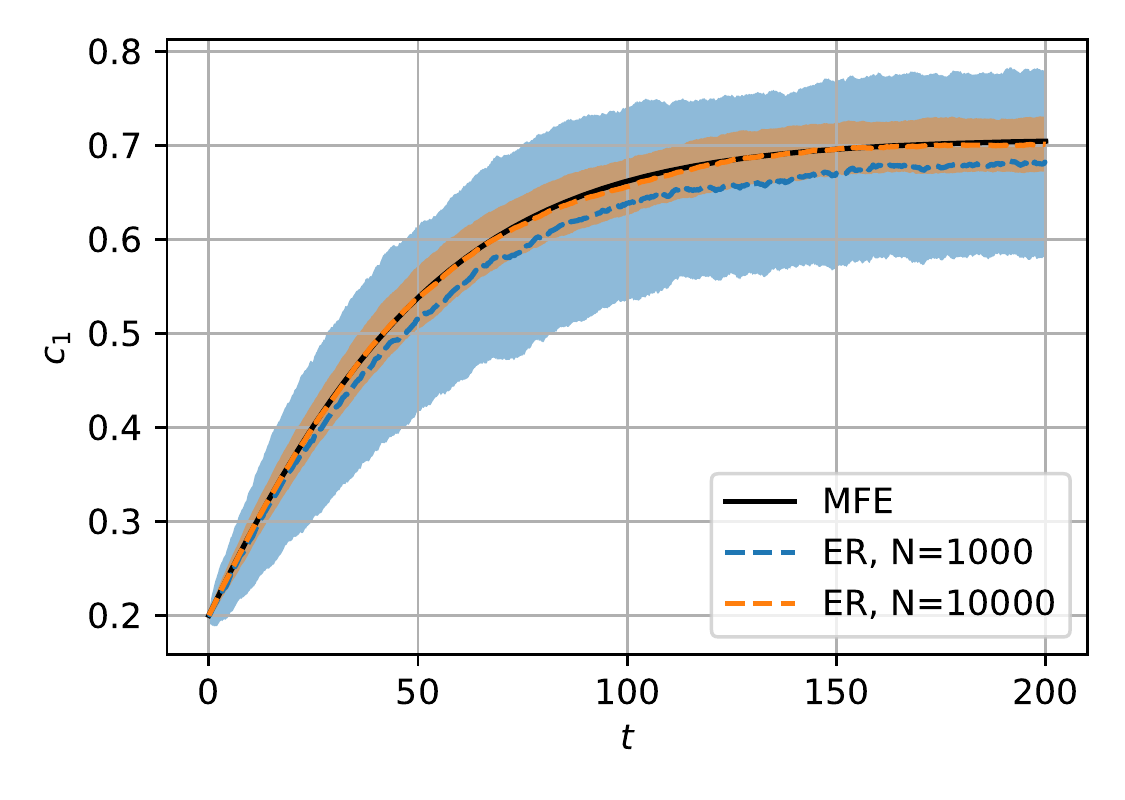}
         \caption{$M=2$ opinions}
         \label{fig:concentration_ER_a}
     \end{subfigure}
     \hfill
     \begin{subfigure}[t]{0.49\textwidth}
         \centering
         \includegraphics[width=\textwidth]{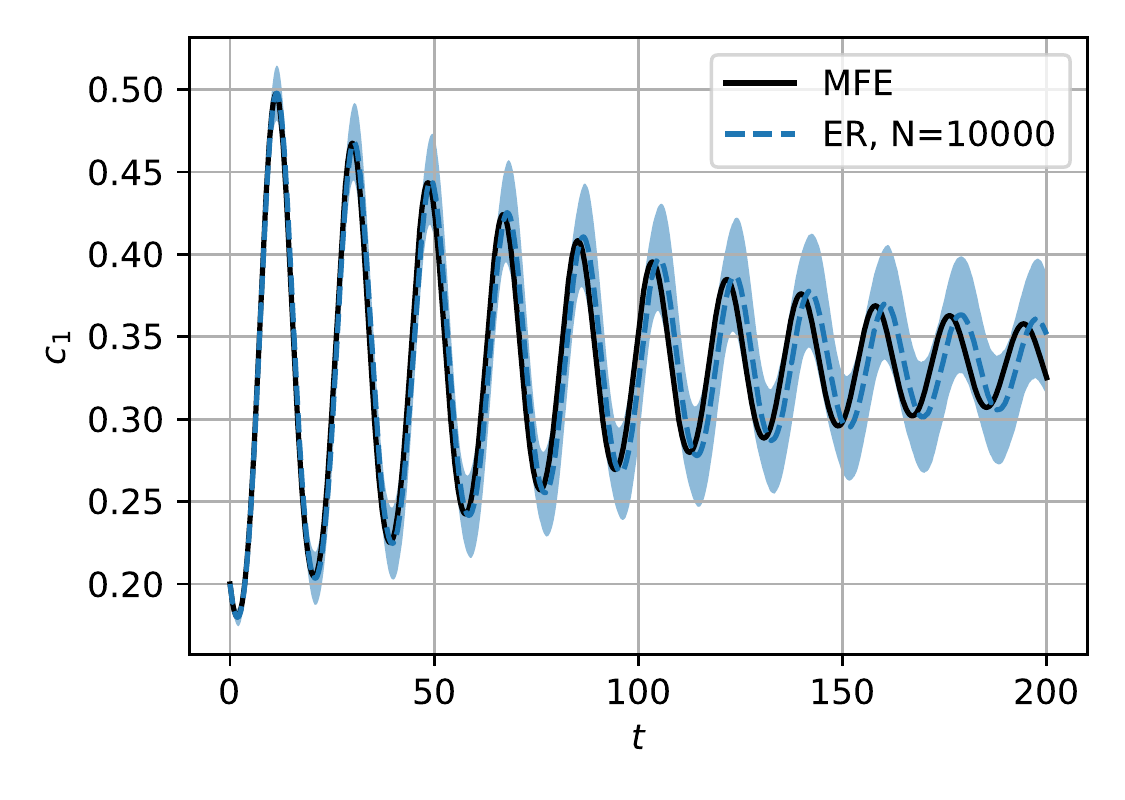}
         \caption{$M=3$ opinions}
         \label{fig:concentration_ER_b}
     \end{subfigure}
        \caption{Mean (dashed line) $\pm$ standard deviation (shaded area) of the CNVM on ER random graphs with edge probability $p=0.01$, estimated from $200$ \review{numerical simulations of the model}, in comparison with the mean-field solution (MFE)~\eqref{eq:mfe_ER_2}. See Example~\ref{example:ER}.}
        \label{fig:concentration_ER}
\end{figure}
\end{example}

\subsection{ER random graphs with heterogeneous population}
\label{subsec:hetero}

Similarly to the previous section, we consider a sequence of Erd\H{o}s--Rényi (ER) random graphs $(\rv{G}_{N,p})_{N\in\mathbb{N}}$.
However, we now assume a heterogeneous population, i.e., there are $K$ distinct classes of agents that differ by their rate constants $r^k_{m,n}$ and $\tilde{r}^k_{m,n}$, $k = 1,\dots,K$, in the CNVM~\eqref{eq:CNVM_rates}. 
For a given graph $G$ an agent $i$ of class $k$ has the rate matrix
\begin{align}
    \big(Q_i^G(x)\big)_{m,n} = r^k_{m, n} \frac{d_{i,n}^G(x)}{d_{i}^G} + \tilde{r}^k_{m,n}\ , \quad m\neq n.
\end{align}
Hence, the collective variable $C_{(m, k)}(x)$ is given by the share of agents that have opinion $m$ and class~$k$.
Note that the quantity of agents in each class and the assignment of agents to the classes are arbitrary, as long as the initial shares $\rv{c}^\ell(0)$ can converge to a constant vector $c_0$ in the large population limit $\ell \to \infty$ (cf.\ Theorem~\ref{thm:main}). This implies that also the shares of each class $k$, i.e., the percentages of agents in each class, have to converge in the large population limit.
Note however that, as discussed in section~\ref{sec:setup}, the class assignment is not allowed to depend on the realization of the random graph.
%In other words, the assignment of the nodes to the classes is arbitrary, but the edges are drawn afterwards and at random.
In other words, for each $\ell$ there is a deterministic assignment of the nodes to the classes, while the edges are drawn afterwards and at random.

Let us again derive the mean-field solution in a heuristic manner.
Consider the propensity functions
\begin{align}
    \alpha^G_{(m,k)\to n}(x) &= \sum_{i: (x_i, s_i) = (m,k)} \big(Q_i^{G}(x)\big)_{m, n} \\
    &= \sum_{i: (x_i, s_i) = (m,k)} \Big( r^k_{m, n}  \frac{d_{i,n}^G(x)}{d_{i}^G} + \tilde{r}^k_{m,n} \Big).
\end{align}
Because of the homogeneous nature of an ER random graph $\rv{G}$, we expect the share of agents of opinion $n$ in the neighborhood of agent $i$, $\smash{ d_{i,n}^{\rv{G}}(x) / d_{i}^{\rv{G}} }$, to be approximately equal to the share of opinion $n$ in the whole system, $c_n := \sum_{k \in [K]} C_{(n, k)}(x)$.
Hence, we have
\begin{align}
    \frac{1}{N} \alpha^{\rv{G}}_{(m,k)\to n}(x) &\approx \frac{1}{N} \sum_{i: (x_i, s_i) = (m,k)} \Big( r^k_{m, n} c_n + \tilde{r}^k_{m,n} \Big) \\
    &= c_{(m,k)} \big(r^k_{m, n} c_n +  \tilde{r}^k_{m,n} \big) =: \tilde{\alpha}_{(m,k)\to n}(c)
\end{align}
which yields the following mean-field ODE when inserted into \eqref{eq:MFE}:
\begin{align} \label{eq:mfe_heterogeneous}
    \frac{\diff}{\diff t} c(t) = \sum_{(m,k) \to n} c_{(m,k)}(t) \Big(r^k_{m, n} \sum_{k^\prime\in[K]} c_{(n, k^\prime)}(t) +  \tilde{r}^k_{m,n} \Big) \big(e_{(n,k)} - e_{(m,k)}\big).
\end{align}

\begin{theorem}
Consider the heterogeneous CNVM as introduced above on the sequence of Erd\H{o}s--Rényi random graphs $(\rv{G}_{N, p})_{N\in\mathbb{N}}$.
Let the edge probability $p=p(N)$ be a function of the number of vertices~$N$.
If $p$ dominates $\log(N)/N$ asymptotically, i.e.,
\begin{align}
    p = \omega\Big(\frac{\log N}{N}\Big) \quad \text{as}\ N \to \infty,
\end{align}
then the dynamics of the collective variables $c = C(x)$, where $c_{(m,k)}$ denotes the share of agents that have opinion $m$ and class $k$, converges to a mean-field limit as $N \to \infty$, \review{in the sense of both Theorem~\ref{thm:main} (annealed result) and Corollary~\ref{cor:quenched} (quenched result)}. The mean-field solution satisfies the ODE~\eqref{eq:mfe_heterogeneous}.
\end{theorem}
\begin{proof}
    The proof is analogous to the proof for a homogeneous population in section~\ref{subsec:ER}.
    We first define
    \begin{align}
        \rv{E}^x_{(m,k)\to n} := \sum_{i: (x_i, s_i) = (m, k)} d_{i, n}^{\rv{G}_{N,p}}(x)
    \end{align}
    as the number of edges between nodes of extended state $(m, k)$ and nodes of opinion~$n$.
    Then, analogously to Lemma~\ref{lemma:ER_counts}, we show by using the Chernoff bound (Lemma~\ref{lemma:chernoff}) that
    \begin{align}\label{Chernoff}
    \mathbb{P}\Big(\Big\lvert \rv{E}^x_{(m,k)\to n} - C_{(m,k)}(x) C_n(x) N^2 p \Big\rvert \geq \varepsilon\Big) \leq 2 \exp\Big(-\frac{\varepsilon^2}{3 N^2 p}\Big),
    \end{align}
    where $C_n(x) := \sum_{k \in [K]} C_{(n, k)}(x)$.
    Now, by inserting the propensity functions to~\eqref{Delta}, we have
    \begin{align}
        \Delta^{\rv{G}}(x) = \max_{(m,k)\to n}\  r^k_{m,n} \Big\lvert \frac{1}{N} \sum_{i : (x_i, s_i) = (m,k)} \frac{d^{\rv{G}}_{i,n}(x)}{d^{\rv{G}}_i} - C_{(m,k)}(x) C_n(x) \Big\rvert.
    \end{align}
    Analogously to the proof of Proposition~\ref{prop:ER_f}, we define the events $\mathcal{A}$ and $\mathcal{B}$, and show by means of~\eqref{Chernoff} that
    \begin{align}
        \mathbb{P}(\mathcal{A}) &\leq \mathbb{P}(\mathcal{A} \cap \mathcal{B}) + \mathbb{P}(\mathcal{B}^C)\\
        &\leq f_\varepsilon(N) := 4 K M^{N+2} \exp\Big(-\frac{1}{12} N^2 p \Big(\frac{\varepsilon}{\hat{r}} - \frac{\varepsilon^2}{\hat{r}^2}\Big)^2 \Big) + 2 N \exp\Big(-N \frac{\varepsilon^2 p}{12 \hat{r}} + \frac{ \varepsilon}{3\hat{r}}\Big).
    \end{align}
    The bounding function $f_\varepsilon$ is identical to the homogeneous case \eqref{f_eps} except for the additional factor $K$, due to the additional maximum over the class $k$ before applying the union bound.
\end{proof}

\begin{figure}
    \centering
    \includegraphics[width=.6\textwidth]{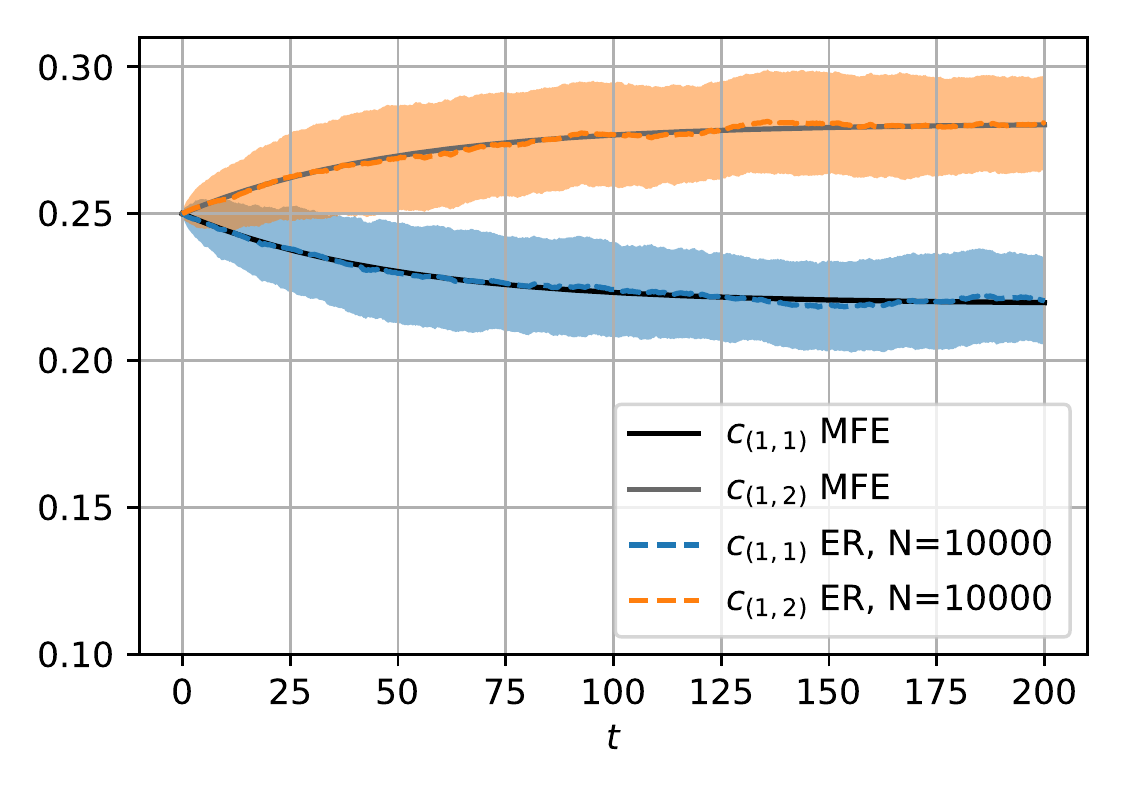}
    \caption{Mean (dashed line) $\pm$ standard deviation (shaded area) of the CNVM ($M=2$ opinions) with heterogeneous population on ER random graphs with $p=0.01$, estimated from $200$ \review{numerical simulations of the model}, in comparison with the mean-field solution~\eqref{eq:mfe_heterogeneous}. The population consists of $K=2$ different classes, and the rates $r^k$ are such that class $1$ slightly prefers opinion 2 and class $2$ prefers opinion $1$, i.e., $r_{1,2}^1 > r_{2,1}^1$ and $r_{1,2}^2 < r_{2,1}^2$.
    Initial conditions are $(c_{(1,1)}, c_{(1,2)}, c_{(2,1)}, c_{(2,2)}) = (0.25, 0.25, 0.25, 0.25)$.}
    \label{fig:concentration_hetero}
\end{figure}

We provide numerical results for an example model in Figure~\ref{fig:concentration_hetero}. 

\subsection{Stochastic block model} \label{subsec:sbm}
In this section we discuss the CNVM (with homogeneous population) on random graphs given by the stochastic block model.
In the stochastic block model the population of agents is split into several clusters (blocks) and there are different edge probabilities for connections inside the clusters and for connections between clusters; see Figure~\ref{fig:illustration} for an example with two clusters.
More precisely, let there be $K$ blocks with sizes $b_1,\dots,b_K \in (0,1] \cap \mathbb{Q}$, such that $\sum_k b_k = 1$.
We consider graphs on $N_\ell = \ell\ \text{LCD}(b_1,\dots,b_k)$ nodes, where LCD refers to the lowest common denominator, and declare that nodes $\{1,\dots, N_\ell b_1\}$ belong to block 1, nodes $\{ N_\ell b_1 + 1,\dots, N_\ell (b_1 + b_2)\}$ to block 2, and so on.
Furthermore, we have a symmetric matrix of probabilities $(p_{k, k^\prime})_{k,k^\prime = 1}^K$, such that $p_{k, k^\prime} \geq 0$ is the probability of an edge between a node in block $k$ and a node in block~$k^\prime$. The edges are then drawn randomly and independently according to these probabilities.
We assume that for all $k \in [K]$ there is at least one $k^\prime \in [K]$ such that~$p_{k, k^\prime} > 0$.

We define the class of node $i$ as $s_i = k$ if node $i$ is located in the $k$-th block.
Hence, the collective variable $C_{(m, k)}(x)$ is given by the share of agents that are located in cluster $k$ and have opinion~$m$.
Let us again derive the mean-field solution in a heuristic manner.
Consider for a given graph $G$ the propensity functions
\begin{align} \label{eq:block_propensity}
    \alpha^{G}_{(m,k)\to n}(x) = \sum_{i: (x_i, s_i) = (m,k)} \Big( r_{m, n}  \frac{d_{i,n}^{G}(x)}{d_{i}^{G}} + \tilde{r}_{m,n} \Big).
\end{align}
(Note that we have a homogeneous population, i.e., every node has equal rate constants $r$ and $\tilde{r}$. It would be straightforward to extend this to the case where nodes have different rate constants $r^k$ and $\tilde{r}^k$ depending on the class $k$, similarly to section~\ref{subsec:hetero}.)
For a random graph $\rv{G}$ generated by the stochastic block model, we expect that the degree $d_i^{\rv{G}}$ of node $i$ in block $k$ is concentrated around
\begin{align}
    d_i^{\rv{G}} \approx N_\ell \sum_{k^\prime \in [K]} b_{k^\prime} p_{k,{k^\prime}} =: N_\ell \bar{p}_k\,.
\end{align}
Thus, we have
\begin{align}
    \alpha^{\rv{G}}_{(m,k)\to n}(x) \approx \frac{r_{m,n}}{N_\ell \bar{p}_k} \sum_{i: (x_i, s_i) = (m,k)} d_{i,n}^{\rv{G}}(x) + \sum_{i: (x_i, s_i) = (m,k)} \tilde{r}_{m, n}\,.
\end{align}
Furthermore, the random variable $\sum_{i: (x_i, s_i) = (m,k)} d_{i,n}^{\rv{G}}(x)$, which counts the number of edges between nodes in cluster $k$ that have opinion $m$ and nodes anywhere that have opinion $n$, is expected to be concentrated around its mean
\begin{align}
    \sum_{i: (x_i, s_i) = (m,k)} d_{i,n}^{\rv{G}}(x) \approx \sum_{{k^\prime} \in [K]} c_{(m,k)} c_{(n,{k^\prime})} N_\ell^2 p_{k, {k^\prime}},
\end{align}
where $c_{(m,k)} = C_{(m,k)}(x)$.
Hence, it follows
\begin{align} 
    \frac{1}{N_\ell} \alpha^{\rv{G}}_{(m,k)\to n}(x) \approx c_{(m,k)}\Big(r_{m, n} \frac{\sum_{{k^\prime} \in [K]} c_{(n,{k^\prime})} p_{k, {k^\prime}}}{\bar{p}_k} + \tilde{r}_{m, n} \Big) =: \tilde{\alpha}_{(m,k)\to n}(c)
    \label{eq:propensity_sbm}
\end{align}
and we obtain, after inserting into \eqref{eq:MFE}, the mean-field ODE
\begin{align} \label{eq:mfe_stochastic_block_model}
    \frac{\diff}{\diff t} c(t) = \sum_{(m,k) \to n} c_{(m,k)}(t)\Big(r_{m, n} \frac{\sum_{{k^\prime} \in [K]} c_{(n,{k^\prime})}(t)\ p_{k, {k^\prime}}}{\bar{p}_k} + \tilde{r}_{m, n} \Big) \big(e_{(n,k)} - e_{(m,k)}\big),
\end{align}
where $\bar{p}_k = \sum_{k^\prime \in [K]} b_{k^\prime} p_{k,{k^\prime}}$.

\begin{theorem} \label{thm:sbm}
Consider the CNVM \eqref{eq:CNVM_rates} on a sequence of stochastic block model random graphs $(\rv{G}_\ell)_{\ell \in \mathbb{N}}$ as defined above.
We introduce a sequence of scaling factors $(\kappa_\ell)_{\ell \in \mathbb{N}}$, $\kappa_\ell \in [0, 1]$, and employ the scaled edge probabilities $\kappa_\ell\, p_{k, {k^\prime}}$ (instead of $p_{k,k'}$) for generating~$\rv{G}_\ell$.
If $\kappa_\ell$ dominates $\log(N_\ell)/N_\ell$ asymptotically, i.e.,
\begin{align}
    \kappa_\ell = \omega\Big(\frac{\log N_\ell}{N_\ell}\Big) \quad \text{as}\ \ell \to \infty, 
\end{align}
then the dynamics of the collective variables $c = C(x)$, where $c_{(m,k)}$ denotes the share of agents that have opinion $m$ and are located in the $k$-th block, converges to a mean-field limit as $\ell \to \infty$, \review{in the sense of both Theorem~\ref{thm:main} (annealed result) and Corollary~\ref{cor:quenched} (quenched result)}. The mean-field solution satisfies the ODE \eqref{eq:mfe_stochastic_block_model}.
\end{theorem}
\begin{proof}
    Again, the proof is analogous to the proof for ER random graphs in section~\ref{subsec:ER}. We provide the detailed proof of this theorem in \cref{sec:appendix_proof_sbm}. The derived bounding function
    \begin{align}
    f_\varepsilon(\ell) = 4 M^{N_\ell + 2} K \exp \Big( -\frac{1}{12} N_\ell^2 \bar{p}^\ell \Big(\frac{\varepsilon}{\hat{r}} -  \frac{\varepsilon^2}{\hat{r}^2}\Big)^2\Big) + 2 N_\ell \exp\Big( -N_\ell \frac{\varepsilon^2 \bar{p}^\ell}{12 \hat{r}} + \frac{\varepsilon}{3 \hat{r}} \Big),
    \label{eq:bounding_sbm}
\end{align}
where $\bar{p}^\ell := \kappa_\ell \min_{k \in [K]} \bar{p}_k$, is identical to the bounding function for ER random graphs (cf.\ Proposition~\ref{prop:ER_f}), except for the additional factor $K$ and the value $\bar{p}^\ell$ instead of~$p$.
\end{proof}

\begin{remark}
In the previous theorem we have considered the case that all edge probabilities $p_{k,k^\prime}$ are scaled using the same factor $\kappa_\ell$.
It is also possible to let each edge probability scale independently, i.e., we define as $p^\ell_{k,k^\prime}$ the edge probabilities used to construct the graph $\rv{G}_\ell$.
For the bounding function $f_\varepsilon(\ell)$ in \eqref{eq:bounding_sbm} to converge to $0$ it is then required that 
\begin{align}
    \bar{p}^\ell := \min_{k\in[K]} \bar{p}_k^\ell := \min_{k\in[K]} \sum_{k^\prime \in [K]} b_{k^\prime} p_{k,k^\prime}^\ell = \omega\Big(\frac{\log N_\ell}{N_\ell}\Big) \, ,
\end{align}
which yields the following condition on the $p^\ell_{k,k^\prime}$ for convergence to a mean-field limit:
\begin{align}
    \forall k \in [K]\ \exists k^\prime \in [K]:\quad p^\ell_{k,k^\prime} = \omega\Big(\frac{\log N_\ell}{N_\ell}\Big) \ \text{as}\ \ell \to \infty.
\end{align}
Moreover, the mean-field equation \eqref{eq:mfe_stochastic_block_model} has to be adapted to this setting: the factor $p_{k,k^\prime} / \bar{p}_k$ in \eqref{eq:mfe_stochastic_block_model} has to be replaced by the limit $\lim_{\ell \to \infty} p^\ell_{k,k^\prime}/\bar{p}_k^\ell$ and hence the edge probabilities $p^\ell_{k,k^\prime}$ may only be chosen in such a way that these limits exist for all $k,k^\prime$.
All in all, this means that for the mean-field limit to hold it is sufficient that every cluster is well connected to at least one other cluster or itself. If a cluster $k$ is only sparsely connected to another cluster $k^\prime$ (or itself), the two are not coupled in the MFE as the factor $p^\ell_{k,k^\prime}/\bar{p}_k^\ell$ becomes $0$ in the limit $\ell \to \infty$.
\end{remark}

\begin{figure}
    \centering
    \includegraphics[width=.6\textwidth]{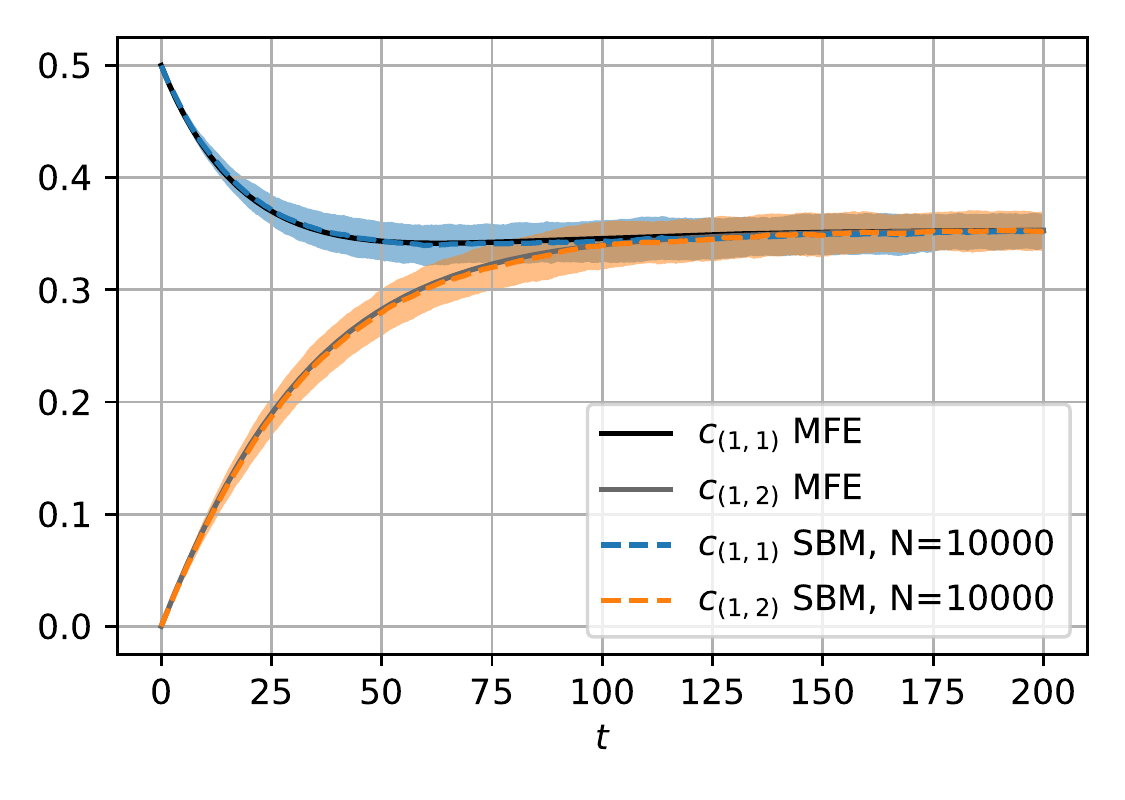}
    \caption{Mean (dotted line) $\pm$ standard deviation (shaded area) of the CNVM on a stochastic block model (SBM) with two equal size blocks and $p_{1,1}= p_{2,2} = 0.01$ and $p_{1,2} = 0.0001$, estimated from $200$ \review{numerical simulations of the model}, in comparison with the mean-field solution \eqref{eq:mfe_stochastic_block_model}. $M=2$ opinions with initial conditions $(c_{(1,1)}, c_{(1,2)}, c_{(2,1)}, c_{(2,2)}) = (0.5, 0, 0, 0.5)$.}
    \label{fig:concentration_sbm}
\end{figure}

We provide numerical results of an example stochastic block model in Figure~\ref{fig:concentration_sbm}.
In the example, there are two equal size blocks and $M=2$ opinions. Initially, every agent in block $1$ has opinion $1$ and every agent in block $2$ has opinion~$2$. Over time the concentrations in both blocks equilibrate.

\subsection{Uniformly random regular graphs} \label{subsec:regular}
In this section we derive mean-field results for the CNVM on uniformly random regular graphs.
A simple graph is called $d$-regular if every node has a degree of exactly~$d$.
We denote by $\rv{G}_{N,d}$ the uniformly random $d$-regular graph on $N$ nodes, i.e., every $d$-regular graph has equal probability and every other graph has probability~$0$.
\review{We again implicitly allow that $d = d(N)$ depends on the size of the graph $N$.}
Similarly to Erd\H{o}s--Rényi random graphs from the previous sections, uniformly random $d$-regular graphs are likely to have a homogeneous edge density, which indicates a mean-field limit with respect to the simple opinion shares ($K=1$), and thus, we employ the same propensity functions and we expect the same mean-field ODE as in section~\ref{subsec:ER}:
\begin{align} \label{eq:MFE_regular}
    \frac{\diff}{\diff t} c(t) = \sum_{m \neq n} c_m(t) \big(r_{m, n} c_n(t) +  \tilde{r}_{m,n} \big) (e_n - e_m).
\end{align}
However, due to the stochastic dependence of edges in $\rv{G}_{N,d}$ (in contrast to ER random graphs) working with random $d$-regular graphs is more intricate, especially in the case of small~$d$.
\review{
In the case of a large degree $d$ on the other hand, the distributions of the random regular graph and the ER random graph with $p=d/N$ become asymptotically identical, which is the subject of the \textit{Sandwich conjecture}~\cite{Kim2004}:
\begin{conjecture}[Sandwich conjecture] \label{sandwich_conjecture}
    If $d = d(N)$ dominates $\log N$ asymptotically, there exist $p_* = (1 - o(1))\ d/N$ and $p^* = (1 + o(1))\ d/N$ as well as ER random graphs $\rv{G}_* \sim \rv{G}_{N, p_*}$ and $\rv{G}^* \sim \rv{G}_{N, p^*}$, such that
    \begin{align}
        \mathbb{P}(\rv{G}_* \subseteq \rv{G}_{N, d} \subseteq \rv{G}^*) = 1 - o(1),
    \end{align}
    where $\subseteq$ denotes inclusion of edges.
\end{conjecture}
\begin{proof}
    To date, the Sandwich conjecture has only been proven for $d \gg \log^4 N / \log^3 \log N$, see~\cite{Gao2020}. It is an open question whether or not the conjecture holds for the missing range $\log N \ll d \ll \log^4 N / \log^3 \log N$.
\end{proof}
Given a degree $d$ large enough such that the above conjecture applies, e.g., $d \approx N^a$ for any fixed $a > 0$, the convergence to the mean-field limit is obtained by our previous Theorem \ref{thm:ER_convergence} for ER random graphs.
\begin{theorem} \label{thm:regular_sandwich}
    Let a sequence of random regular graphs $(\rv{G}_{N_\ell, d})_{\ell \in \mathbb{N}}$ be given such that $d = d(N_\ell) \gg \log^4 N_\ell / \log^3 \log N_\ell$. Then the dynamics of the opinion shares in the CNVM converges to a mean-field limit as $\ell \to \infty$, in the sense of both Theorem~\ref{thm:main} (annealed result) and Corollary~\ref{cor:quenched} (quenched result). The associated mean-field ODE is given in~\eqref{eq:MFE_regular}.
\end{theorem}
}

As stated before, the case of random regular graphs with small degree $d$ is substantially more difficult.
It is often easier to deal with the configuration model instead, which we introduce below, and then transfer the results back to the regular simple graph setting.

\paragraph{Configuration model}
The configuration model \cite[section 11.1]{Frieze2015} enables us to define random multigraphs with arbitrary degree distributions. However, we restrict ourselves to $d$-regular multigraphs.
Unlike simple graphs, multigraphs allow self-loops, i.e., edges from a node to itself, and also multiple edges between two nodes. For our purposes, we will need to discard multigraphs that are not simple from the sampling.

The configuration model is constructed as follows.
We choose $N$ and $d$ such that $Nd$ is even and define the set $W = [Nd]$.
Consider the partition $W = W_1 \cup \dots \cup W_N$, where $W_i := \{(i-1)d + 1, \dots, i d \}$.
Moreover, we define the map
\begin{equation}
    \label{eq:phiassign}
    \varphi:[Nd]\to[N],\ \varphi(e) = i :\Leftrightarrow e \in W_i.
\end{equation}
In this setup, the elements of $W_i$ refer to half-edges attached to the $i$-th node, and joining two half-edges $e,h \in W$ denotes forming an edge between the nodes $\varphi(e)$ and $\varphi(h)$ in the multigraph.
More precisely, we call a partition $F$ of $W$ into $\eta := Nd / 2$ pairs a \textit{configuration} and denote the multigraph $G=(V,E)$ that is induced by the configuration $F$ by
\begin{align}
    \gamma(F) := (V, E) = \Big([N], \{(\varphi(e), \varphi(h)) \mid (e,h) \in F\}\Big).
\end{align}
Let the random variable $\rv{F}$ denote a uniformly random configuration, i.e., every possible configuration $F$ is equally likely.
Then one can show that for any two $d$-regular simple graphs $G_1$ and $G_2$ we have
\begin{align}
    \mathbb{P}(\gamma(\rv{F}) = G_1) = \mathbb{P}(\gamma(\rv{F}) = G_2),
\end{align}
and hence, $\rv{G}_{N,d}$ can be obtained by conditioning the configuration model $\gamma(\rv{F})$ on the set of simple graphs~\cite[Corollary~11.2]{Frieze2015}.

\begin{figure}
\centering
\includegraphics[width=.9\textwidth]{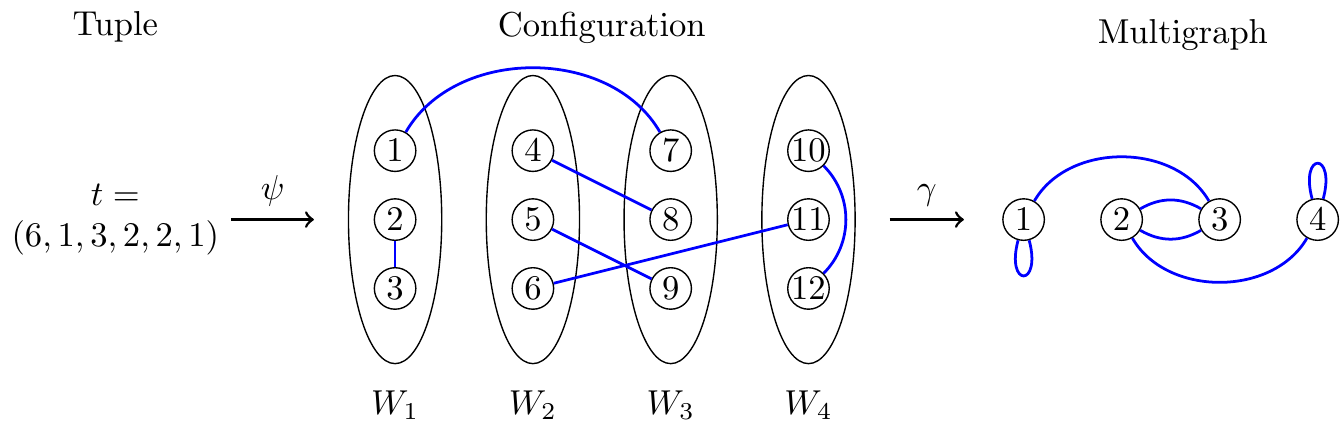}
\caption{For $d=3$ and $N=4$ and an example tuple $t \in \Pi$, we show the induced configuration $\psi(t)$ and multigraph~$\gamma(\psi(t))$.}
\label{fig:configuration_model}
\end{figure}

It remains to find a simple way to sample configurations $F$ uniformly at random.
Let
\begin{align}
    \Pi := \Big\{(t_1, \dots, t_\eta) \mid\forall r\in [\eta]:\ t_r \in \{1,\dots, (Nd -2r + 1)\}\Big\}.
\end{align}
A tuple $t \in \Pi$ uniquely induces a configuration $F = \psi(t)$, where the map $\psi$ is defined via the following procedure.
Let $U^1 := W$ and $U^{r+1} := U^r \setminus \{u_0^r, u_{t_r}^r\}$, $r=1,\dots,\eta$, where $u_i^r$ is the $(i+1)\text{-th}$ smallest element of~$U^r$.
Then $F = \psi(t) := \{(u_0^r, u_{t_r}^r) \mid r=1,\dots,\eta\}$. 
In words, we start with $W$ and define an edge by connecting the nodes associated to the first and $t_1$-th element of~$W$. Then we remove this pair of elements from $W$ and continue this procedure on the remaining set with the first and $t_2$-th element, and so on.
We provide an example in Figure~\ref{fig:configuration_model}.
Note that, for the random variable $\rv{t} \in \Pi$ that is uniformly distributed, i.e., each component $t_r \in \{1,\dots, (Nd -2r + 1)\}$ is picked uniformly at random and independently from the others, we have
\begin{align}
    \psi(\rv{t}) \overset{d}{=} \rv{F}.
\end{align}

In order to verify the conditions of Theorem \ref{thm:main} for uniformly random $d$-regular graphs, the following concentration result is useful:
\begin{lemma} \label{lemma:concentration_configurationmodel}
Let $x \in [M]^N$ and fix two distinct opinions $m,n \in [M]$. Assume that the state $x$ is ordered, such that the $m$-opinion nodes are first, the $n$-opinion nodes come after that, and then the rest.
Define the random variable $g(\rv{t})$ as the number of edges between nodes of opinion $m$ and nodes of opinion $n$ in the induced multigraph $\gamma(\psi(\rv{t}))$, with respect to~$x$. Then it follows
\begin{align}
        \mathbb{P}\Bigg( \bigg|g(\rv{t}) -  \frac{c_m c_n N^2 d^2}{Nd - 1}\bigg| \geq \varepsilon \Bigg) \leq 2 \exp\Big(-\frac{\varepsilon^2}{4 Nd}\Big),
\end{align}
where $c_m$ denotes the share of opinion $m$ in the state $x$.
\end{lemma}
\begin{proof}
We assume that there is at least one node with opinion $m$ and at least one node with opinion $n$ in $x$; otherwise the lemma is trivially true.
Consider two tuples $t,t^\prime \in \Pi$ that only differ in one coordinate $\ldif$, i.e., $t = (t_1, \dots, t_\eta)$, $t^\prime = (t_1,\dots,t_{\ldif - 1}, t_\ldif^\prime, t_{\ldif + 1}, \dots, t_\eta)$.
Let $b_m \in W$ denote the maximal element such that $x_{\varphi(b_m)} = m$ (recall~\eqref{eq:phiassign}), and define $b_n$ analogously.
Note that, due to the ordering of $x$, we have~$b_m < b_n$.
The values $b_m$ and $b_n$ act as important boundaries because the edges counting towards $g(t)$ have to cross $b_m$ but must not cross $b_n$.
(An edge $(s,e)$ with $s,e \in W$ is defined to be crossing the boundary $b$ if $s \leq b < e$.)
In Lemma \ref{lemma:A_boundary_crossing}, we show that the number of edges crossing any boundary can vary by at most $2$ between $t$ and~$t^\prime$.
Hence, it follows that~$|g(t) - g(t^\prime)| \leq 4$, because there are the two boundaries $b_m,b_n$ to consider for~$g(t)$.
As the random vector $\rv{t} = (\rv{t}_1, \dots, \rv{t}_\eta)$ has independent components, we are able to apply McDiarmid's inequality~\cite{McDiarmid1989}:
\begin{align}
        \mathbb{P}\Big( \Big| g(\rv{t}) - \mathbb{E}[g(\rv{t})] \Big| \geq \varepsilon \Big) \leq 2 \exp\Big(-\frac{2 \varepsilon^2}{16 \eta}\Big) = 2 \exp\Big(-\frac{ \varepsilon^2}{4 Nd}\Big).
    \end{align}
Finally, note that there are $c_m N d$ half-edges attached to nodes of opinion $m$, and each of these has a $c_n N d / (Nd - 1)$ chance to get matched with a half-edge of a node with opinion~$n$.
Hence, we have
\begin{align}
    \mathbb{E}[g(\rv{t})] = c_m N d \frac{c_n N d}{Nd - 1},
\end{align}
which completes the proof.
\end{proof}

Now we can derive the bounding function $f_\varepsilon$ from the conditions of Theorem~\ref{thm:main}. We consider a sequence of uniformly random regular graphs $\rv{G}_\ell = \rv{G}_{N_\ell, d}$, $\ell \in \mathbb{N}$. Note that for a fixed degree $d\in\mathbb{N}$ not all graph sizes $N\in\mathbb{N}$ are possible, hence the sequence $(N_\ell)_{\ell}$ is necessary.
\begin{proposition} \label{prop:regular_bound}
For all $\varepsilon > 0$ there exists a function $f_\varepsilon: \mathbb{N} \to \mathbb{R}_{\geq 0}$ such that
\begin{align}
    \forall \ell \in \mathbb{N}:\ \mathbb{P}\Big(\max_{x\in [M]^{N_\ell}} \Delta^{\rv{G}_{\ell}}(x) \geq \varepsilon \Big) \leq f_\varepsilon(\ell),
\end{align}
where
\begin{align}
    f_\varepsilon(\ell) = (2+o(1)) M^{N_\ell+2} \exp\Big(d^2 - N_\ell d \frac{\varepsilon^2}{4 \hat{r}^2} + \frac{\varepsilon}{4 \hat{r}}\Big) \quad \text{as }\ell\to\infty.
\end{align}
\end{proposition}
\begin{proof}
Let $\ell \in \mathbb{N}$ and denote~$\rv{G} := \rv{G}_\ell$. By inserting the propensity functions (cf.\ equations~\eqref{eq:cnvm_propensity} and \eqref{eq:cnvm_reduced_propensity}), we have
\begin{align}
        \Delta^{\rv{G}}(x) = \max_{m\neq n}\  r_{m,n} \Big\lvert \frac{1}{N_\ell} \sum_{i : x_{i} = m} \frac{d^{\rv{G}}_{i,n}(x)}{d^{\rv{G}}_i} - C_m(x) C_n(x) \Big\rvert.
    \end{align}
Fix two opinions $m\neq n$ and let $x \in [M]^{N_\ell}$ be ordered as in Lemma~\ref{lemma:concentration_configurationmodel}. For simpler notation, we write $c_m := C_m(x)$. As all realizations of $\rv{G}$ are $d$-regular, it follows
\begin{align}
    \frac{1}{N_\ell} \sum_{i : x_{i} = m} \frac{d^{\rv{G}}_{i,n}(x)}{d^{\rv{G}}_i} = \frac{1}{N_\ell d} \sum_{i : x_{i} = m} d^{\rv{G}}_{i,n}(x) =: \frac{1}{N_\ell d} \rv{E}_{m,n} \, ,
\end{align}
where $\rv{E}_{m,n}$ denotes the number of edges between nodes of opinion $m$ and nodes of opinion $n$ in $\rv{G}_{N_\ell, d}$, with respect to $x$.
Consider also the number of edges between nodes of opinion $m$ and nodes of opinion $n$ in the configuration model, which we denote by $g(\rv{t})$ as in Lemma~\ref{lemma:concentration_configurationmodel}.
We can relate these two quantities as follows, provided $d < N_\ell^{1/7}$~\cite[Thm 11.3]{Frieze2015}:
\begin{align}
        \mathbb{P}\Big(\Big|\rv{E}_{m,n} - \frac{c_m c_n N_\ell^2 d^2}{N_\ell d - 1}\Big| \geq \varepsilon\Big)  &\leq   (1+o(1)) e^{\lambda(\lambda+1)} \mathbb{P}\Big(\Big|g(\rv{t}) - \frac{c_m c_n N_\ell^2 d^2}{N_\ell d - 1}\Big| \geq \varepsilon\Big) \\
        &\leq (1+o(1)) e^{\lambda(\lambda+1)} 2 \exp\Big(-\frac{\varepsilon^2}{4 N_\ell d}\Big) \label{eq:concentration_regular}
\end{align}
where $\lambda := \frac{d - 1}{2}$. With the notation $\hat{r} := \max_{m\neq n} r_{m,n}$, we have
\begin{align}
    &\mathbb{P}\Big( r_{m,n} \Big\lvert \frac{1}{N_\ell} \sum_{i : x_{i} = m} \frac{d^{\rv{G}}_{i,n}(x)}{d^{\rv{G}}_i} - c_{m}c_{n} \Big\rvert \geq \varepsilon\Big) \\
    &\leq \mathbb{P}\Big( \Big\lvert \frac{1}{N_\ell} \sum_{i : x_{i} = m} \frac{d^{\rv{G}}_{i,n}(x)}{d^{\rv{G}}_i} - \frac{c_m c_n N_\ell d}{N_\ell d - 1} \Big\rvert + \Big\lvert \frac{c_m c_n N_\ell d}{N_\ell d - 1} - c_{m}c_{n} \Big\rvert \geq \frac{\varepsilon}{\hat{r}} \Big)\\
    &= \mathbb{P}\Big( \frac{1}{N_\ell d} \Big\lvert  \rv{E}_{m,n} - \frac{c_m c_n N_\ell^2 d^2}{N_\ell d - 1}  \Big\rvert + \frac{c_m c_n}{N_\ell d - 1} \geq \frac{\varepsilon}{\hat{r}} \Big)\\
    &= \mathbb{P}\Big(\Big\lvert \rv{E}_{m,n} -\frac{c_m c_n N_\ell^2 d^2}{N_\ell d - 1}  \Big\rvert \geq N_\ell d \big(\frac{\varepsilon}{\hat{r}} - \frac{c_m c_n}{N_\ell d - 1}\big) \Big)\\
    &\overset{\eqref{eq:concentration_regular}}{\leq} (1+o(1)) e^{\lambda(\lambda+1)} 2 \exp\Big(-\frac{\big(N_\ell d (\frac{\varepsilon}{\hat{r}} - \frac{c_m c_n}{N_\ell d - 1})\big)^2}{4 N_\ell d}\Big)\\
    &\leq (1+o(1)) e^{\lambda(\lambda+1)} 2 \exp\Big( - N_\ell d \frac{\varepsilon^2}{4 \hat{r}^2} + \frac{\varepsilon}{4 \hat{r}} \Big). \label{eq:regular_bound}
\end{align}
Recall that we have assumed an ordered state $x$. However, due to the indifference of the random regular graph with respect to the specific node numbering, a certain regular graph is just as likely as the same graph but with permuted node labels. Using this property, it follows that the bound \eqref{eq:regular_bound} holds for general states $x$. We give a more detailed explanation in \cref{sec:appendix_isomorph}.
Finally, we apply the union bound
\begin{align}
    &\mathbb{P}\Big(\max_{x\in [M]^{N_\ell}} \Delta^{\rv{G}}(x) \geq \varepsilon \Big)\\
    &= \mathbb{P}\Big(\max_{x \in [M]^{N_\ell}} \max_{m \neq n} r_{m,n} \Big\lvert \frac{1}{N_\ell} \sum_{i : x_{i} = m} \frac{d^{\rv{G}}_{i,n}(x)}{d^{\rv{G}}_i} - C_m(x) C_n(x) \Big\rvert \geq \varepsilon\Big)\\
    &\overset{\eqref{eq:regular_bound}}{\leq} (1 + o(1)) M^{N_\ell} M(M-1)  e^{\lambda(\lambda+1)} 2 \exp\Big( - N_\ell d \frac{\varepsilon^2}{4\hat{r}^2} + \frac{\varepsilon}{4\hat{r}} \Big) \\
    &\leq (2+o(1)) M^{N_\ell+2} \exp\Big(d^2 - N_\ell d \frac{\varepsilon^2}{4 \hat{r}^2} + \frac{\varepsilon}{4 \hat{r}}\Big)
\end{align}
to conclude the proof.
\end{proof}

As in the previous sections, we can summarize our findings in the following
\begin{theorem} \label{thm:regular}
Consider the CNVM \eqref{eq:CNVM_rates} on a sequence of uniformly random regular graphs $\rv{G}_\ell = \rv{G}_{N_\ell, d_\ell}$, $\ell \in \mathbb{N}$.
If $d_\ell \to \infty$ as $\ell \to \infty$, but slower than $N_\ell^{1/7}$, i.e.
\begin{align}
    d_\ell = \omega(1) \cap o(N_\ell^{1/7}),
\end{align}
then the dynamics of the opinion shares concentrates around a mean-field limit as $\ell \to \infty$, \review{in the sense of both Theorem~\ref{thm:main} (annealed result) and Corollary~\ref{cor:quenched} (quenched result)}. The mean-field solution satisfies the ODE \eqref{eq:MFE_regular}.
\end{theorem}
\begin{proof}
We need to guarantee that the bounding function derived in Proposition \ref{prop:regular_bound}
\begin{align}
    f_\varepsilon(\ell) = (2+o(1)) M^{N_\ell+2} \exp\Big(d_\ell^2 - N_\ell d_\ell \frac{\varepsilon^2}{4 \hat{r}^2} + \frac{\varepsilon}{4 \hat{r}}\Big)
\end{align}
converges to $0$ as $\ell \to \infty$. Hence, the exponent (after removing constants) $N_\ell + d_\ell^2 - N_\ell d_\ell \varepsilon^2$ has to converge to $-\infty$ for all $\varepsilon > 0$. Using $N_\ell \to \infty$ and $d_\ell \geq 1$, we require
\begin{align}
    \lim_{\ell \to \infty} \frac{N_\ell + d_\ell^2}{N_\ell d_\ell \varepsilon^2} = \lim_{\ell \to \infty} \frac{1}{d_\ell \varepsilon^2} + \frac{d_\ell}{N_\ell \varepsilon^2} \overset{!}{=} 0.
\end{align}
The left term goes to $0$ for $d_\ell = \omega(1)$, and the right term for $d_\ell = o(N_\ell)$. Additionally, we used \cite[Thm 11.3]{Frieze2015} in the previous proposition, which requires that  $d_\ell = o(N_\ell^{1/7})$.

\review{Moreover, neglecting constants, $f_\varepsilon(\ell)$ is bounded by $\exp(N_\ell + d_\ell^2 - N_\ell d_\ell) \ll N_\ell^{-2}$, from which the condition $\sum_\ell f_\varepsilon(\ell) < \infty$ of Corollary~\ref{cor:quenched} follows.}
\end{proof}

\review{
Combining our insights for random regular graphs with small degree and with large degree yields the following corollary.
\begin{corollary} \label{cor:regular}
    Let a sequence of random regular graphs $(\rv{G}_{N_\ell, d_\ell})_{\ell \in \mathbb{N}}$ satisfy 
    \begin{align}
        \lim_{\ell \to \infty} d_\ell = \infty.
    \end{align}
    Then the dynamics of the opinion shares in the CNVM converges to a mean-field limit as $\ell \to \infty$, in the sense of both Theorem~\ref{thm:main} (annealed result) and Corollary~\ref{cor:quenched} (quenched result). The associated mean-field ODE is given in~\eqref{eq:MFE_regular}.
\end{corollary}
\begin{proof}
     We have shown the convergence to the mean-field limit for $1 \ll d \ll N^{1/7}$ in Theorem~\ref{thm:regular} and for $d \gg \log^4 N / \log^3 \log N$ in Theorem~\ref{thm:regular_sandwich}.
     Due to the fact that $N^{1/7}\gg \log^4 N / \log^3 \log N$, the convergence follows for all $d \gg 1$.
\end{proof}
}

\begin{remark}
    In section \ref{subsec:ER} we showed convergence to the MFE for Erd\H{o}s--Rényi random graphs if the average node degree grows faster than $\log(N)$, whereas our result for $d$-regular graphs (Theorem~\ref{thm:regular}) only requires unboundedness of the degree $d$, i.e., it can grow arbitrarily slowly. Intuitively, the regularity of the graphs allows for the mean-field limit to hold under more sparsity.
\end{remark}

\begin{figure}
     \centering
     \begin{subfigure}[t]{0.49\textwidth}
         \centering
         \includegraphics[width=\textwidth]{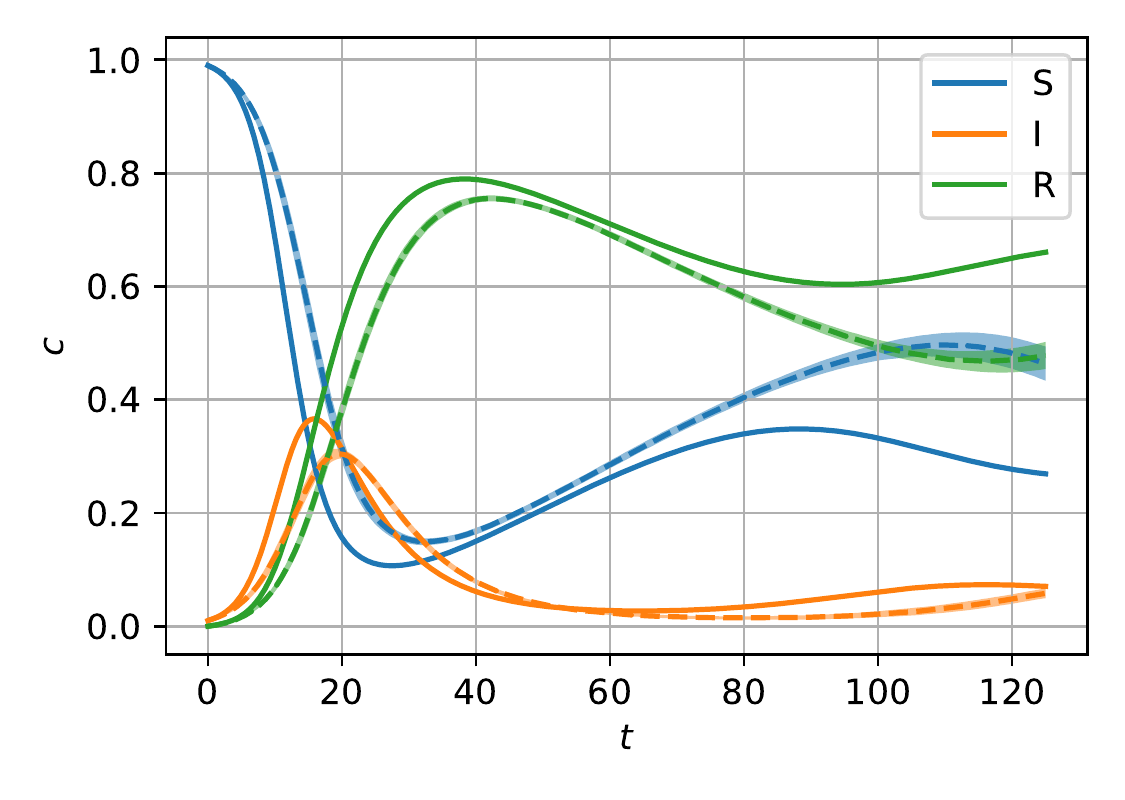}
         \caption{Degree $d = 10$.}
     \end{subfigure}
     \hfill
     \begin{subfigure}[t]{0.49\textwidth}
         \centering
         \includegraphics[width=\textwidth]{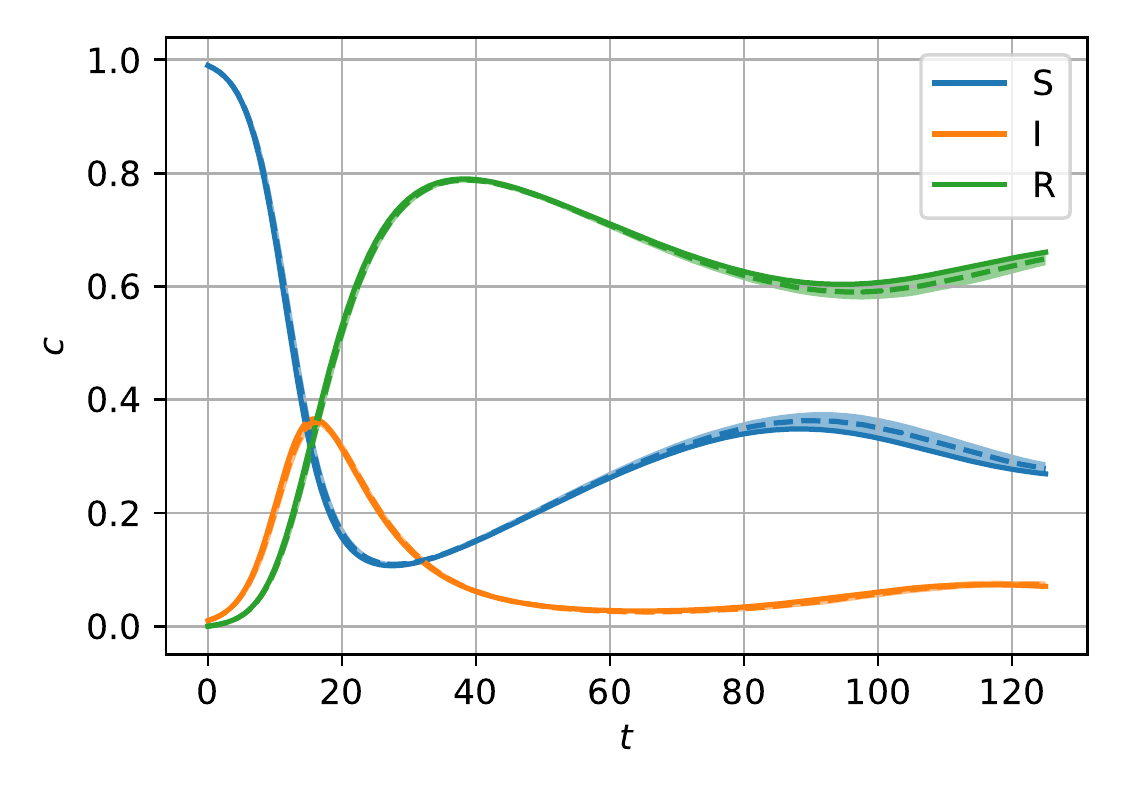}
         \caption{Degree $d = 100$.}
     \end{subfigure}
        \caption{Mean (dotted line) $\pm$ standard deviation (shaded area) of the SIR model (cf. Example \ref{example:SIR}) on uniformly random regular graphs with $N=10000$ nodes, estimated from $100$ realizations, in comparison to the mean-field solution \eqref{eq:MFE_regular} (solid line). Initial shares $(c_S, c_I, c_R) = (0.99, 0.01, 0)$.}
        \label{fig:concentration_regular}
\end{figure}

We provide a numerical example for uniformly random regular graphs in Figure~\ref{fig:concentration_regular}.
We choose the rates $r$ and $\tilde{r}$ to replicate the SIR model discussed in Example~\ref{example:SIR} and observe a steep wave of infections followed by a smaller second wave.
The figure illustrates how the discrepancy between model realizations and mean-field solution decreases when we increase the degree $d$, as indicated by Theorem~\ref{thm:regular}.
For $d = 10$, the approximation quality of the mean-field solution is poor, even though the number of agents $N=10000$ is quite large. (Increasing $N$ reduces the variance of the model realizations, but does not necessarily move the mean closer to the mean-field solution.)
For $d = 100$, on the other hand, the mean-field limit is a reasonable approximation.
Hence, in order to achieve a certain approximation quality it is crucial that both $N$ and $d$ are large enough.

\section{Conclusion} \label{sec:conclusion}
In this work we have derived conditions under which Markovian discrete-state systems on networks (as described in section \ref{sec:setup}) converge to a mean-field limit. More precisely, we have provided an ordinary differential equation (ODE) such that the shares of nodes of certain classes behave like the solution of this ODE in the large population limit, if the conditions are fulfilled.
Moreover, we have applied these results to the well-known voter model on Erd\H{o}s--Rényi random graphs, the stochastic block model, and uniformly random regular graphs, specifying the convergence conditions for each of these graph types. As for Erd\H{o}s--Rényi graphs we also have shown how to incorporate a heterogeneous population with different classes of agents. 
Numerical examples have validated the derived mean-field solutions.

While we have provided verifiable conditions that guarantee convergence to a mean-field limit, the question of finding optimal collective variables for a given problem is still open. In many (rather simple) examples, like the ones we have discussed in this paper, it is obvious which collective variables to choose, i.e., on which subsets of nodes the concentrations of states should be measured. In order to deal with more intricate problems, methods to construct viable collective variables could be investigated in future work. \review{One of the first of such methods was proposed in~\cite{Luecke2023}, where also scale-free networks are studied. Indeed, such network topologies that better resemble the social network structures found in real-life data~\cite{Mislove2007} are of particular interest and could be subject of future theoretical studies}.

Furthermore, fluctuations around the mean-field solution could be analyzed in order to investigate the convergence to the mean-field limit more thoroughly. For medium-sized populations, that do not allow a good approximation via the mean-field limit, these fluctuations could be taken into account to construct an approximation via a stochastic differential equation, e.g., by adding a stochastic correction term to the mean-field limit~\cite{Niemann2021}.

Moreover, the framework introduced in this paper could be extended to more general models in future works. For instance, a changing network structure that is coupled to the dynamics of the nodes' states, e.g., as in the adaptive voter-model \cite{PhysRevE.74.056108}, could be investigated.

\section*{Acknowledgments}
This work has been supported by Deutsche Forschungsgemeinschaft (DFG) under
Germany’s Excellence Strategy via the Berlin Mathematics Research Center MATH+ (EXC2046/ project ID: 390685689).
The authors thank Michael Anastos for valuable discussions about random graphs\review{, and they thank the anonymous reviewer for excellent suggestions that allowed to improve the paper.}

\newpage
% \begin{appendices}
\appendix
\setcounter{equation}{0}
\renewcommand{\theequation}{A.\arabic{equation}}
\section{Auxiliary lemmas}
\begin{lemma} \label{lemma:conv_prob}
Let $(\Omega, \mathcal{F}, \mathbb{P})$ be a probability space and $\rv{z}\super{\ell}: \Omega \times \R \to [0, 1]$, $\ell \in \mathbb{N}$, a sequence of stochastic processes with Lebesgue-measurable realizations, i.e., $\rv{z}\super{\ell}(\omega, \cdot)$ is measurable for all $\omega$. We denote the random variable describing the process at time $t \in \R$ as $\rv{z}\super{\ell}(t)$.
Assume that for all $\varepsilon > 0$ there exists a function $f_\varepsilon:\mathbb{N} \to \R$ such that
\begin{enumerate}[a)]
    \item $\forall \ell \in \mathbb{N} \ \forall s \in \R: \mathbb{P}(\rv{z}\super{\ell}(s) \geq \varepsilon) \leq f_\varepsilon(\ell)$
    \item $\lim_{\ell \to \infty} f_\varepsilon(\ell) = 0.$
\end{enumerate}
Then it follows that
\begin{enumerate}[1)]
    \item $\forall s \in \R: \lim_{\ell \to \infty} \mathbb{E}[\rv{z}\super{\ell}(s)] = 0$
    \item $\forall t\in\R_{\geq 0}: \lim_{\ell \to \infty} \mathbb{E}[\int_0^t \rv{z}\super{\ell}(s) ds] = 0$
    \item $\forall t\in\R_{\geq 0}: \int_0^t \rv{z}\super{\ell}(s) \diff s \toprob 0$ as $\ell \to \infty$.
\end{enumerate}
\end{lemma}
\begin{proof}
Let $\varepsilon > 0$ and define for every $\ell \in \mathbb{N}$ the stochastic process $\hat{\rv{z}}\super{\ell}(s)$ by
\begin{align}
    \hat{\rv{z}}\super{\ell}(s) := \begin{cases}
    \varepsilon, & \rv{z}\super{\ell}(s) < \varepsilon\\
    1, & \text{else.}
    \end{cases}
\end{align}
Thus, we have $\rv{z}\super{\ell}(s) \leq \hat{\rv{z}}\super{\ell}(s)$ and
\begin{align}
    \mathbb{E}[\rv{z}\super{\ell}(s)] &\leq \mathbb{E}[\hat{\rv{z}}\super{\ell}(s)] = \varepsilon \mathbb{P}(\rv{z}\super{\ell}(s) < \varepsilon) + \mathbb{P}(\rv{z}\super{\ell}(s) \geq \varepsilon)\\
    &\leq \varepsilon + f_\varepsilon(\ell).
\end{align}
By assumption b) this yields $\lim_{\ell \to \infty} \mathbb{E}[\rv{z}\super{\ell}(s)] \leq \varepsilon$.
As $\varepsilon > 0$ was arbitrary, statement 1) follows.\\
In order to prove statement~2), we first use Tonelli's theorem to interchange integral and expected value, i.e.
\begin{align}
    \mathbb{E}\Big[\int_0^t \rv{z}\super{\ell}(s) ds\Big] = \int_0^t \mathbb{E}[\rv{z}\super{\ell}(s)] ds,
\end{align}
since the integrand is non-negative.
Given $\mathbb{E}[\rv{z}\super{\ell}(s)] \leq 1$, we can apply the dominated convergence Theorem, which yields
\begin{align}
    \lim_{\ell \to \infty} \int_0^t \mathbb{E}[\rv{z}\super{\ell}(s)] ds = \int_0^t \lim_{\ell \to \infty} \mathbb{E}[\rv{z}\super{\ell}(s)] ds = 0
\end{align}
by statement~1). Statement 3) follows directly from 2), as convergence in $L^1$ is stronger than convergence in probability.
\end{proof}

\begin{lemma}[Chernoff bound] \label{lemma:chernoff}
Let $\rv{X}_1, \dots, \rv{X}_n$ be independent random variables with values in $\{0, 1\}$ and denote $\rv{X} := \sum_i \rv{X}_i$. Then for all $\varepsilon > 0$
\begin{align}
    \mathbb{P}\Big(\big\lvert \rv{X} - \mathbb{E}[\rv{X}] \big\rvert \geq \varepsilon\Big) \leq 2 \exp\Big(-\frac{\varepsilon^2}{3\ \mathbb{E}[\rv{X}]}\Big).
\end{align}
\begin{proof}
    See for example \cite[Theorem A.14]{Arora2009}.
\end{proof}

\end{lemma}

\begin{lemma} \label{lemma:A_ER_degrees}
Let $\rv{G}_{N,p}$ denote the ER random graph and the random variable $\rv{d}_i := d_i^{\rv{G}_{N,p}}$ the degree of the $i$-th node.
Then for all $\varepsilon > 0$ and all $i \in [N]$ we have
\begin{align}
    \mathbb{P}(|\rv{d}_i - Np| \geq \varepsilon Np) \leq 2 \exp\Big( -\frac{\varepsilon^2 N p}{3} + \frac{2 \varepsilon}{3}\Big).  
\end{align}
\end{lemma}
\begin{proof}
    In $\rv{G}_{N,p}$ the degree $\rv{d}_i$ of each vertex is Binomial distributed with $N-1$ trials and success probability $p$.
    Hence, we have $\mathbb{E}[\rv{d}_i] = (N-1)p$ and using the Chernoff bound (Lemma \ref{lemma:chernoff}) it follows that
    \begin{align}
        \mathbb{P}\Big(|\rv{d}_i - (N-1)p| \geq \varepsilon \Big) \leq 2 \exp\Big(-\frac{\varepsilon^2}{3 (N-1) p}\Big).
        \label{eq:degree_concentration}
    \end{align}
    Thus, we have
    \begin{align}
        \mathbb{P}\Big(|\rv{d}_i - Np| \geq \varepsilon N p \Big) &\leq \mathbb{P}\Big(|\rv{d}_i - (N-1) p| + p \geq \varepsilon N p \Big)\\
        &\overset{\eqref{eq:degree_concentration}}{\leq} 2 \exp\Big(-\frac{(\varepsilon N p - p)^2}{3 (N-1) p}\Big)\\
        &\leq 2 \exp\Big(-\frac{(\varepsilon N p - p)^2}{3Np}\Big)\\
        &\leq 2 \exp\Big(-\frac{\varepsilon^2 N p}{3} + \frac{2 \varepsilon}{3}\Big).
    \end{align}
\end{proof}

\begin{lemma} \label{lemma:A_boundary_crossing}
Let $1\leq b \leq Nd$ and for a tuple $t \in \Pi$ let $h(t) := |\{(s, e) \in \psi(t) \mid s \leq b < e\}|$ denote the number of edges that cross the boundary $b$.
Let $t, t^\prime \in \Pi$ only differ in one coordinate $\ldif$, i.e., $t = (t_1, \dots, t_\eta)$, $t^\prime = (t_1,\dots,t_{\ldif - 1}, t_\ldif^\prime, t_{\ldif + 1}, \dots, t_\eta)$.
Then it follows that
\begin{align}
    \lvert h(t) - h(t^\prime) \rvert \leq 2.
\end{align}
\end{lemma}
\begin{proof}
    Let $U^r$ and $u_i^r$ be as introduced in section~\ref{subsec:regular}. Define
\begin{equation}
\label{eq:delta_def}
    \delta^r := \begin{cases}
        1, & \text{if } u_0^r \leq b < u_{t_r}^r\\
        0, & \text{else}
    \end{cases}
\end{equation}
and note that $h(t) = \sum_{r=1}^\eta \delta^r$.
Moreover, let $i^r := \max \{i \mid u_i^r \leq b\}$ be the index of the largest element in $U^r$ that is not larger than~$b$. (We set $\max \emptyset := -1$.)
The following relations between $\delta^r$ and $i^r$ hold:
\begin{align}
0 \leq i^r < t_r \quad &\Leftrightarrow \quad \delta^r = 1 \label{eq:delta_i_1}\\
    i^r \geq t_r \text{ or } i^r = -1 \quad &\Leftrightarrow \quad \delta^r = 0 \label{eq:delta_i_2}\\
    \delta^r = 1 \quad &\Rightarrow \quad i^{r+1} = i^r - 1 \label{eq:delta_i_3}\\
    \delta^r = 0 \quad &\Rightarrow \quad i^{r+1} =
    \begin{cases}
    -1, & \text{if } i^r = -1\\
    i^r - 2, & \text{else.} \label{eq:delta_i_4}
    \end{cases}
\end{align}
Relation \eqref{eq:delta_i_3} holds because $\delta^r = 1$ implies that exactly one element that is smaller or equal to $b$ is removed from $U^r$, and hence $i^{r+1}$ is one less than $i^{r}$.
Relation \eqref{eq:delta_i_4} holds because $\delta^r = 0$ implies that either two elements that are smaller or equal to $b$ are removed from $U^r$, which yields a reduction of $i^{r+1}$ by 2 compared to $i^{r}$, or both removed elements are larger than $b$, which is the case if $i^r = -1$.\\
We denote the respective analogous objects for $t^\prime$ as $(U^\prime)^r$, $(u^\prime)_i^r$, $ (\delta^\prime)^r $, and $(i^\prime)^r$.
W.l.o.g., we assume that $t_\ldif < t_\ldif^\prime$. \\
Clearly, for $r \leq \ldif$ we have $U^r = (U^\prime)^r$ and $i^r = (i^\prime)^r$, and for $r < \ldif$ we have $\delta^{r} = (\delta^\prime)^{r}$.
Consider the case that it also holds $\delta^{\ldif} = (\delta^\prime)^{\ldif}$. Due to equations \eqref{eq:delta_i_3} and \eqref{eq:delta_i_4}, this yields $i^{\ldif+1} = (i^\prime)^{\ldif+1}$, and as $t_r = t^\prime_r$ for $r>\ldif$, we have $\delta^{\ldif + 1} = (\delta^\prime)^{\ldif + 1}$ via equations \eqref{eq:delta_i_1} and \eqref{eq:delta_i_2}. By iteration, we have $\delta^r = (\delta^\prime)^{r}$ for all $r>\ell$ and hence $h(t) = h(t^\prime)$.\\
We now consider the case that $\delta^{\ldif} \neq (\delta^\prime)^{\ldif}$, i.e., by $t_\ldif < t_\ldif^\prime$ and~\eqref{eq:delta_def}, $\delta^{\ldif}=0$ and $(\delta^\prime)^{\ldif}=1$. By \eqref{eq:delta_i_2} this implies $i^\ldif = (i^\prime)^{\ldif} \neq -1$.
Also, by equations \eqref{eq:delta_i_3} and \eqref{eq:delta_i_4} we have $(i^\prime)^{\ldif + 1}= i^{\ldif + 1} + 1$.
Depending on the value of $(i^\prime)^{\ldif + 1}$, one of the following three outcomes occurs:
\begin{enumerate}
    \item
    $(i^\prime)^{\ldif + 1} = 0$: It follows from~\eqref{eq:delta_i_1} that $(\delta^\prime)^{\ldif + 1} = 1$, and because $i^{\ldif + 1} = -1$ we have $\delta^{\ldif + 1} = 0$.
    Moreover, in the next step from~\eqref{eq:delta_i_3} we have $(i^\prime)^{\ldif + 2} = -1 = i^{\ldif + 2}$ and thus, $(i^\prime)^{r} = i^{r}$ also for all subsequent steps $r > \ldif + 2$.
    As a result, for all $r > (\ldif + 1)$ we have $\delta^r = (\delta^\prime)^{r}=0$.
    All in all, this yields $h(t^\prime) = h(t) + 2$.
    \item
    $(i^\prime)^{\ldif + 1} = t_{\ldif +1}$: It follows from~\eqref{eq:delta_i_2} that $(\delta^\prime)^{\ldif + 1} = 0$, and because $i^{\ldif + 1} = t_{\ldif +1} - 1$ we have $\delta^{\ldif + 1} = 1$.
    Hence, by equations \eqref{eq:delta_i_3} and \eqref{eq:delta_i_4} we have $(i^\prime)^{\ldif + 2} = i^{\ldif + 2}$ and thus, by iteration also $(i^\prime)^{r} = i^{r}$ for all subsequent steps $r > \ldif + 2$.
    As a result, for all $r > (\ldif + 1)$ we have $\delta^r = (\delta^\prime)^{r}$.
    All in all, this yields $h(t^\prime) = h(t) + 1 - 1 = h(t)$.
    \item
    Else: For all other values of $(i^\prime)^{\ldif + 1}$, we have $\delta^{\ldif + 1} = (\delta^\prime)^{\ldif + 1}$. Thus, in the next step we still have the relation $(i^\prime)^{\ldif + 2}= i^{\ldif + 2} + 1$ due to \eqref{eq:delta_i_3} if $\delta^{\ldif + 1} = 1$, or due to \eqref{eq:delta_i_4} if $\delta^{\ldif + 1} = 0$.
    By iteration, we get into either case 1.\ or 2.\ in one of the subsequent steps $r > (\ldif + 1)$, i.e., either $(i^\prime)^{r} = 0$ or $(i^\prime)^{r} = t_{r}$, which yields $h(t^\prime) = h(t) + 2$, or $h(t^\prime) = h(t)$ respectively.
\end{enumerate}
\end{proof}

\section{Proof of Theorem \ref{thm:main}} \label{sec:appendix_proof_main}
Let $\big(\rv{P}_{(m,k)\to n}(t)\big)_{m,k,n}$ denote independent unit-rate Poisson processes. Then we can write (cf. \cite[section 1.2]{Winkelmann2020})
\begin{align}
    \rv{c}^\ell(t) = \rv{c}^\ell(0) + \sum_{(m,k)\to n} \rv{P}_{(m,k)\to n} \Big( \int_{0}^{t} \alpha_{(m,k)\to n}^{\rv{G}_\ell}(\rv{x}^\ell(s))\, \diff s\Big) \frac{v_{(m,k)\to n}}{N_\ell}.
\end{align}
We further define the centered Poisson processes $\rv{\tilde{P}}_{(m,k)\to n}(t) := \rv{P}_{(m,k)\to n}(t) - t$, which leads to
\begin{align}
    \rv{c}^\ell(t) = \rv{c}^\ell(0) + \underbrace{
    \sum_{(m,k)\to n} \rv{\tilde{P}}_{(m,k)\to n} \Big( \int_{0}^{t} \alpha_{(m,k)\to n}^{\rv{G}_\ell}(\rv{x}^\ell(s))\, \diff s\Big) \frac{v_{(m,k)\to n}}{N_\ell}
    }_{=: \rv{\delta}_{\ell}(t)} + \int_0^t \rv{F}_{\ell}(\rv{x}^\ell(s)) \diff s,
    \label{eq:rv_c}
\end{align}
where
\begin{align}
    \rv{F}_{\ell}(x) := \sum_{(m,k)\to n} \alpha_{(m,k)\to n}^{\rv{G}_\ell}(x) \frac{v_{(m,k)\to n}}{N_\ell}.
\end{align}
Note that due to the assumption \eqref{eq:Q_bound_B} of transition rates bounded by $B>0$, we have
\begin{align}
    \alpha^{\rv{G}_\ell}_{(m,k)\to n}(x) = \sum_{i: (x_i, s_i) = (m,k)} \big(Q_i^{\rv{G}_\ell}(x)\big)_{m, n} \leq N_\ell B
\end{align}
and thus
\begin{align}
    \hat{\rv{\delta}}_\ell(t) &:= \sup_{0 \leq s \leq t}\  \lVert \rv{\delta}_{\ell}(s) \rVert \\
    &\leq \sum_{(m,k)\to n} \sup_{0 \leq s \leq t} \Big\lvert \frac{1}{N_\ell} \rv{\tilde{P}}_{(m,k)\to n} (s N_\ell B) \Big\rvert \lVert v_{(m,k)\to n} \rVert. \label{eq:hat_delta}
\end{align}
By the law of large numbers, one can show that (see for example \cite[Theorem~1.2]{Anderson2011})
\begin{align}
    \sup_{0 \leq s \leq t} \big\lvert \frac{1}{N_\ell} \rv{\tilde{P}}_{(m,k)\to n} (s N_\ell B) \big\rvert \toprob 0 \quad \text{as } \ell \to \infty
\end{align}
 and hence
\begin{align}
     \forall t:\ \hat{\rv{\delta}}_\ell(t) \underset{\ell \to \infty}{\toprob} 0.
     \label{eq:delta_conv}
\end{align}
Furthermore, we have that
\begin{align}
    &\Big\lVert \int_0^t \rv{F}_{\ell}(\rv{x}^\ell(s)) - F(\rv{c}^\ell(s)) ds \Big\rVert\\
    &\leq \int_0^t \underbrace{
    \sum_{(m,k)\to n} \Big\lvert \frac{1}{N_\ell} \alpha_{(m,k)\to n}^{\rv{G}_\ell}(\rv{x}^\ell(s)) - \tilde{\alpha}_{(m,k)\to n}(\rv{c}^\ell(s)) \Big\rvert \lVert v_{(m,k)\to n} \rVert
    }_{=: \rv{z}^\ell(s)} \diff s =: \rv{\tilde{\delta}}_\ell(t) \label{eq:tilde_delta_def}
\end{align}
and
\begin{align}
    \rv{z}^\ell(s) \leq \sum_{(m,k)\to n} \Delta^{\rv{G}_\ell}\big(\rv{x}^\ell(s)\big) \bar{v} = MK(M-1)\ \Delta^{\rv{G}_\ell}\big(\rv{x}^\ell(s)\big) \bar{v},
\end{align}
where $\bar{v} := \max_{(m,k)\to n} \lVert v_{(m,k)\to n} \rVert$. Let $\varepsilon > 0$ and define $\tilde{\varepsilon} := \frac{\varepsilon}{MK(M-1)\bar{v}}$. Then it follows that
\begin{align}
    \Pb(\rv{z}^\ell(s) \geq \varepsilon) &\leq \Pb\Big(\Delta^{\rv{G}_\ell}\big(\rv{x}^\ell(s)\big) \geq \tilde{\varepsilon}\Big)\\
    &\leq \Pb \Big(\max_{x \in [M]^{N_\ell}} \Delta^{\rv{G}_\ell}(x) \geq \tilde{\varepsilon}\Big)\\
    &\leq f_{\tilde{\varepsilon}}(\ell) \overset{\ell \to \infty}{\longrightarrow} 0.
\end{align}
Hence, from Lemma \ref{lemma:conv_prob} it follows that
\begin{align}
    \rv{\tilde{\delta}}_\ell(t) = \int_0^t \rv{z}^\ell(s) \diff s \underset{\ell  \to \infty}{\toprob}0.
    \label{eq:delta_tilde_conv}
\end{align}
Now, writing $c(t) = c(0) + \int_0^t F(c(s)) \diff s$ and $\rv{c}^\ell(t)$ as in \eqref{eq:rv_c}, we obtain
\begin{align}
    \lVert \rv{c}^\ell(t) - c(t) \rVert &= \Big\lVert \rv{c}^\ell(0) - c(0) +  \rv{\delta}_\ell(t) + \int_0^t \rv{F}_{\ell}(\rv{x}^\ell(s)) - F(c(s)) \diff s \Big\rVert\\
    &\leq \lVert \rv{c}^\ell(0) - c(0) \rVert + \lVert \rv{\delta}_\ell(t) \rVert + \Big\lVert \int_0^t \rv{F}_{\ell}(\rv{x}^\ell(s)) - F(\rv{c}^\ell(s)) \diff s \Big\rVert\\
    &\qquad\qquad\qquad\qquad\qquad\quad+ \Big\lVert \int_0^t F(\rv{c}^\ell(s)) - F(c(s)) \diff s \Big\rVert \nonumber\\
    & \leq \lVert \rv{c}^\ell(0) - c(0) \rVert + \hat{\rv{\delta}}_\ell(t) + \rv{\tilde{\delta}}_\ell(t) + L \int_0^t \rVert \rv{c}^\ell(s) - c(s) \rVert \diff s,
\end{align}
where $L$ denotes the Lipschitz constant of $F$. ($F$ is Lipschitz continuous because we have assumed that all $\tilde{\alpha}_{(m,k)\to n}$ are.)
Note that $\hat{\rv{\delta}}_\ell(t)$ and $\rv{\tilde{\delta}}_\ell(t)$ are monotonically increasing in $t$. Thus, by the Gronwall lemma we obtain
\begin{align}
    \lVert \rv{c}^\ell(t) - c(t) \rVert \leq \Big(\lVert \rv{c}^\ell(0) - c(0) \rVert + \hat{\rv{\delta}}_\ell(t) + \rv{\tilde{\delta}}_\ell(t) \Big) \exp(L t) \underset{\ell \to \infty}{\toprob 0} \label{eq:approx_quality}
\end{align}
due to \eqref{eq:delta_conv} and \eqref{eq:delta_tilde_conv}.
Because of the monotonicity of the above bound, the theorem follows.

\review{
\begin{remark} \label{remark_appendix_convergence}
    We show that the rate of convergence of $\hat{\rv{\delta}}_\ell(t)$ to $0$ as $\ell \to \infty$ is $\sqrt{N_\ell}^{-1}$, in the sense that $\mathbb{E}[\hat{\rv{\delta}}_\ell(t)] =  O(\sqrt{N_\ell}^{-1})$.
    Neglecting constant factors, we can conclude from \eqref{eq:hat_delta} that $\mathbb{E}[\hat{\rv{\delta}}_\ell(t)]$ is bounded by 
    \begin{align}
        \mathbb{E}\left[\sup_{0 \leq s \leq t} \Big\lvert \frac{1}{N_\ell} \rv{\tilde{P}} (s N_\ell B) \Big\rvert\right],
    \end{align}
    where $\rv{\tilde{P}}$ is a centered Poisson process and $B>0$ is the bound of the transition rates, see~\eqref{eq:Q_bound_B}.
    One can show \cite[Lemma~1.3]{Anderson2011} that the centered Poisson process is approximated well by Brownian motion $\rv{W}$, i.e., for all $\ell \in \mathbb{N}$,
    \begin{align}
        \rv{\Gamma} := \sup_{t \geq 0} \frac{\lvert \rv{\tilde{P}}(t N_\ell B) - \rv{W}(t N_\ell B) \rvert}{\log(\max(2,t N_\ell B))} < \infty \quad \text{a.s.,} \qquad \mathbb{E}[\rv{\Gamma}] < \infty,
    \end{align}
    from which we have (assuming $\ell$ large enough so that $tN_\ell B \geq 2$) 
    \begin{align}
        \sup_{0\leq s \leq t} \Big\lvert \rv{\tilde{P}} (s N_\ell B) - \sqrt{N_\ell}\; \rv{W}(sB) \Big\rvert \leq \rv{\Gamma} \log(t N_\ell B).
    \end{align}
    This implies
    \begin{align}
        \sup_{0 \leq s \leq t} \Big\lvert \frac{1}{N_\ell} \rv{\tilde{P}} (s N_\ell B) \Big\rvert 
        \leq \rv{\Gamma} \frac{\log(t N_\ell B)}{N_\ell} + \sup_{0 \leq s \leq t} \Big\lvert \frac{1}{\sqrt{N_\ell}} \rv{W}(sB) \Big\rvert
    \end{align}
    and
    \begin{align}
        \mathbb{E}\left[\sup_{0 \leq s \leq t} \Big\lvert \frac{1}{N_\ell} \rv{\tilde{P}} (s N_\ell B) \Big\rvert\right] \leq \mathbb{E}[\rv{\Gamma}] \frac{\log(t N_\ell B)}{N_\ell} + \frac{1}{\sqrt{N_\ell}} \mathbb{E}\left[\sup_{0 \leq s \leq t} \lvert \rv{W}(sB) \rvert\right],
    \end{align}
    from which the claim $\mathbb{E}[\hat{\rv{\delta}}_\ell(t)] =  O(\sqrt{N_\ell}^{-1})$ follows.
\end{remark}
}

\section{Proof of Proposition \ref{prop:ER_f}} \label{sec:appendix_proof_ER_f}
    We fix any $N \in \mathbb{N}$ and denote $\rv{G} := \rv{G}_{N,p}$.
    By inserting the propensity functions \eqref{eq:cnvm_propensity} and \eqref{eq:cnvm_reduced_propensity}, it follows that
    \begin{align}
        \Delta^{\rv{G}}(x) = \max_{m\neq n}\  r_{m,n} \Big\lvert \frac{1}{N} \sum_{i : x_{i} = m} \frac{d^{\rv{G}}_{i,n}(x)}{d^{\rv{G}}_i} - C_m(x) C_n(x) \Big\rvert.
    \end{align}
    Let $\delta \in (0,1)$ and define the events
    \begin{align}
        \mathcal{A} &:= \Big\{ \max_{x\in [M]^{N}} \Delta^{\rv{G}_N}(x) \geq \varepsilon \Big\} = \Big\{ \max_{x\in [M]^{N}} \max_{m \neq n}\ r_{m,n} \Big\lvert \frac{1}{N} \sum_{i : x_{i} = m} \frac{d^{\rv{G}}_{i,n}(x)}{d^{\rv{G}}_i} - C_m(x) C_n(x) \Big\rvert \geq \varepsilon \Big\},\\
        \mathcal{B} &:= \Big\{ \forall i: (1-\delta)Np \leq d^{\rv{G}}_i \leq (1+\delta)Np\Big\}.
    \end{align}
    From equation~\eqref{eq:ER_degrees} and the union bound, it follows that
    \begin{align} \label{eq:ER_sandwich}
        \mathbb{P}(\mathcal{B}^C) \leq 2 N \exp\Big(\frac{1}{3} (-\delta^2 N p + 2 \delta)\Big),
    \end{align}
    where $\mathcal{B}^C$ denotes the complement of $\mathcal{B}$. We will now show that $\mathbb{P}(\mathcal{A} \cap \mathcal{B})$ is small, from which the Proposition follows when combined with \eqref{eq:ER_sandwich}.
    We define $\rv{E}^x_{m,n} := \sum_{i:x_i=m} d_{i,n}^{\rv{G}}(x)$ as we did in Lemma \ref{lemma:ER_counts}. For any fixed state $x \in [M]^N$ and any opinions $m \neq n$, we have
    \begin{align}
        &\mathbb{P}\Big( r_{m,n} \Big\lvert \frac{1}{N} \sum_{i : x_{i} = m} \frac{d^{\rv{G}}_{i,n}(x)}{(1\pm \delta)Np} - C_m(x) C_n(x) \Big\rvert \geq \varepsilon \Big)\\
        &= \mathbb{P}\Big( \frac{1}{(1\pm \delta) N^2 p} r_{m,n} \big\lvert \rv{E}^x_{m,n} - C_m(x) C_n(x) (1\pm \delta) N^2 p \big\rvert \geq \varepsilon \Big)\\
        & \leq \mathbb{P}\Big( \big\lvert \rv{E}^x_{m,n} - C_m(x) C_n(x) N^2 p \big\rvert + \underbrace{C_m(x) C_n(x)}_{\leq 1} \delta N^2 p \geq r_{m,n}^{-1}  (1\pm \delta) N^2 p \varepsilon \Big)\\
        & \leq \mathbb{P}\Big(\big\lvert \rv{E}^x_{m,n} - C_m(x) C_n(x) N^2 p \big\rvert \geq r_{m,n}^{-1} (1\pm \delta) N^2 p \varepsilon - \delta N^2 p  \Big) \\
        & \leq \mathbb{P}\Big( \big\lvert \rv{E}^x_{m,n} - C_m(x) C_n(x) N^2 p \big\rvert \geq N^2 p (\hat{r}^{-1} \varepsilon - \hat{r}^{-1} \varepsilon \delta - \delta) \Big).
    \end{align}
    Hence, this also holds after applying the maximum:
    \begin{align}
        &\mathbb{P}\Big(\max_{x\in [M]^N} \max_{m\neq n}\ r_{m,n} \Big\lvert \frac{1}{N} \sum_{i : x_{i} = m} \frac{d^{\rv{G}}_{i,n}(x)}{(1\pm \delta)Np} - C_m(x) C_n(x) \Big\rvert \geq \varepsilon \Big)\\
        &\leq \mathbb{P}\Big(\max_{x\in [M]^N} \max_{m\neq n} \big\lvert \rv{E}^x_{m,n} - C_m(x) C_n(x) N^2 p \big\rvert \geq N^2 p (\hat{r}^{-1} \varepsilon - \hat{r}^{-1} \varepsilon \delta - \delta) \Big).
    \end{align}
    In order to ensure that $\hat{r}^{-1} \varepsilon - \hat{r}^{-1} \varepsilon \delta - \delta > 0$ for the given $\varepsilon \in (0,\hat{r})$, we choose $\delta = \hat{r}^{-1} \varepsilon / 2$, i.e., $(\hat{r}^{-1} \varepsilon - \hat{r}^{-1} \varepsilon \delta - \delta) = (\hat{r}^{-1} \varepsilon - \hat{r}^{-2} \varepsilon^2)/2$.
    Applying the union bound and Lemma \ref{lemma:ER_counts} yields
    \begin{align}
        &\mathbb{P}\Big(\max_{x\in [M]^{N}} \max_{m \neq n}\ r_{m,n} \Big\lvert \frac{1}{N} \sum_{i : x_{i} = m} \frac{d^{\rv{G}}_{i,n}(x)}{(1\pm \delta)Np} - C_m(x) C_n(x) \Big\rvert \geq \varepsilon \Big) \\
        &\leq 2 M^N M(M-1) \exp\Big(-\frac{(N^2 p (\hat{r}^{-1} \varepsilon - \hat{r}^{-2} \varepsilon^2)/2))^2}{3 N^2 p}  \Big)  \\
        &\leq 2 M^{N+2} \exp\Big(-\frac{1}{12} N^2 p \Big(\frac{\varepsilon}{\hat{r}} - \frac{\varepsilon^2}{\hat{r}^2}\Big)^2 \Big). \label{eq:ER_estim_bounds}
    \end{align}
    Moreover, we have
    \begin{align}
        &\mathbb{P}(\mathcal{A} \cap \mathcal{B})\\
        &= \mathbb{P}\Big( \max_{x\in [M]^{N}} \max_{m \neq n}\ r_{m,n} \Big\lvert \frac{1}{N} \sum_{i : x_{i} = m} \frac{d^{\rv{G}}_{i,n}(x)}{d^{\rv{G}}_i} - C_m(x) C_n(x) \Big\rvert \geq \varepsilon \\
        & \qquad \qquad \text{and}\ \forall i: (1-\delta)Np \leq d^{\rv{G}}_i \leq (1+\delta)Np \Big) \nonumber\\
        &\leq \mathbb{P}\Big(\max_{x\in [M]^{N}} \max_{m \neq n} \max_{\xi_1,\dots,\xi_N \in [-1,1]} r_{m,n} \Big\lvert \frac{1}{N} \sum_{i : x_{i} = m} \frac{d^{\rv{G}}_{i,n}(x)}{(1 + \xi_i \delta) N p} - C_m(x) C_n(x) \Big\rvert \geq \varepsilon \Big).
    \end{align}
    Since the right term, $C_m(x) C_n(x)$, is independent of $i$, the maximum is reached by either making the left term, $\sum_{i : x_{i} = m} {d^{\rv{G}}_{i,n}(x)} / {(1 + \xi_i \delta) N p}$, as large as possible or as small as possible, i.e., either all $\xi_i = 1$ or all $\xi_i = -1$.
    Therefore, we have, again using the union bound
    \begin{align}
        \mathbb{P}(\mathcal{A} \cap \mathcal{B}) &\leq \mathbb{P}\Big( \max_{\xi \in \{-1,1\}} \max_{x\in [M]^{N}} \max_{m \neq n} r_{m,n} \Big\lvert \frac{1}{N} \sum_{i : x_{i} = m} \frac{d^{\rv{G}}_{i,n}(x)}{(1 + \xi \delta) N p} - C_m(x) C_n(x) \Big\rvert \geq \varepsilon \Big)  \\
       &\overset{\eqref{eq:ER_estim_bounds}}{\leq} 4 M^{N+2} \exp\Big(-\frac{1}{12} N^2 p \Big(\frac{\varepsilon}{\hat{r}} - \frac{\varepsilon^2}{\hat{r}^2}\Big)^2 \Big). \label{eq:ER_conditional}
    \end{align}
    Finally, it follows
    \begin{align}
        \mathbb{P}(\mathcal{A}) &\leq \mathbb{P}(\mathcal{A} \cap \mathcal{B}) + \mathbb{P}(\mathcal{B}^C)\\
        &\hspace{-2.6ex}\overset{\eqref{eq:ER_sandwich}, \eqref{eq:ER_conditional}}{\leq} 4 M^{N+2} \exp\Big(-\frac{1}{12} N^2 p \Big(\frac{\varepsilon}{\hat{r}} - \frac{\varepsilon^2}{\hat{r}^2}\Big)^2 \Big) + 2 N \exp\Big(-N \frac{\varepsilon^2 p}{12 \hat{r}} + \frac{ \varepsilon}{3\hat{r}}\Big).
    \end{align}

\section{Proof of Theorem~\ref{thm:sbm} (convergence for the stochastic block model)} \label{sec:appendix_proof_sbm}
In this section we verify the conditions of Theorem \ref{thm:main} for the continuous-time noisy voter model on stochastic block model random graphs.
The proof is analogous to the proof for ER random graphs in section \ref{subsec:ER}.
For simplicity of notation, we consider the edge probabilities $p_{k,{k^\prime}}$ without the scaling factor~$\kappa_\ell$.
We begin with an analogous version of Lemma \ref{lemma:ER_counts}.
\begin{lemma} \label{lemma:sbm_counts}
Given a fixed state $x \in [M]^N$ and the stochastic block model random graph $\rv{G} = \rv{G}_{N_\ell}$, we define
\begin{align}
    \rv{E}^x_{(m,k)\to n} := \sum_{i: (x_i, s_i) = (m, k)} d_{i, n}^{\rv{G}}(x)
\end{align}
    as the number of edges between nodes of extended state $(m, k)$ and nodes of opinion $n$.
    Then we have
\begin{align}
    \mathbb{E}[\rv{E}^x_{(m,k)\to n}] = \sum_{{k^\prime} \in [K]} C_{(m,k)}(x) C_{(n,{k^\prime})}(x) N^2 p_{k,{k^\prime}} =: \mu
\end{align}
and, using the notation $\bar{p}_k := \sum_{{k^\prime} \in [K]} b_{k^\prime} p_{k,{k^\prime}}$, we have
\begin{align}
    \mathbb{P}\Big( \Big\lvert \rv{E}^x_{(m,k)\to n} - \mu \Big\rvert \geq \varepsilon\Big) \leq 2 \exp\Big(-\frac{ \varepsilon^2}{3 N^2 \bar{p}_k }\Big).
\end{align}
\end{lemma}
\begin{proof}
    The number of edges between a node with extended state $(m,k)$ and a node with extended state $(n,{k^\prime})$ is binomial distributed with $C_{(m,k)}(x) C_{(n,{k^\prime})}(x) N^2$ trials and success probability $p_{k,{k^\prime}}$, i.e.,
    \begin{align}
        \rv{E}^x_{(m,k)\to n} \sim \sum_{{k^\prime} \in [K]} \text{Bin}\big(C_{(m,k)}(x) C_{(n,{k^\prime})}(x) N^2, p_{k,{k^\prime}}\big).
    \end{align}
    From the Chernoff bound (Lemma \ref{lemma:chernoff}), we have
    \begin{align}
    \mathbb{P}\Big( \Big\lvert \rv{E}^x_{(m,k)\to n} - \mu \Big\rvert\geq \varepsilon\Big) \leq 2 \exp\Big(-\frac{\varepsilon^2}{3 \mu}\Big) \leq 2 \exp\Big(-\frac{ \varepsilon^2}{3 N^2 \bar{p}_k }\Big),
    \end{align}
    where the last inequality is due to $C_{(m,k)}(x) C_{(n,{k^\prime})}(x) \leq b_{k^\prime}$.
\end{proof}

Moreover, we show that the node degrees are concentrated, analogously to Lemma \ref{lemma:A_ER_degrees}.

\begin{lemma} \label{lemma:sbm_degrees}
Let node $i \in [N]$ be in cluster $k$ and let $\rv{d}_i := d_i^{\rv{G}}$ denote the degree of node $i$ in the stochastic block model.
Then for all $\varepsilon > 0 $ we have
\begin{align}
    \mathbb{P}\Big(\lvert \rv{d}_i - N_\ell \bar{p}_k) \rvert \geq \varepsilon N_\ell \bar{p}_k \Big) \leq 2 \exp\Big(- N_\ell \frac{\varepsilon^2 \bar{p}}{3} + \frac{2 \varepsilon}{3} \Big),
\end{align}
where $\bar{p} := \min_{k\in [K]} \bar{p}_k$.
\end{lemma}
\begin{proof}
    Note that $\rv{d_i}$ is the sum of independent binomial random variables
    \begin{align}
        \rv{d}_i \sim \sum_{{k^\prime} \in [K]\setminus \{k\}} \text{Bin}(N_\ell b_{k^\prime}, p_{k, {k^\prime}}) + \text{Bin}(N_\ell b_k - 1, p_{k, k}).
    \end{align}
    Using the abbreviation $\mu := \mathbb{E}[\rv{d}_i]$, we have $N_\ell \bar{p}_k = \mu + p_{k, k}$ and 
    \begin{align}
        \mathbb{P}\Big(\lvert \rv{d}_i - N_\ell \bar{p}_k) \rvert \geq \varepsilon N_\ell \bar{p}_k \Big)
        &\leq \mathbb{P}\Big(\lvert \rv{d}_i - \mu) \rvert + p_{k,k} \geq \varepsilon N_\ell \bar{p}_k \Big)\\
        &\leq 2 \exp\Big(- \frac{(\varepsilon N_\ell \bar{p}_k - p_{k,k})^2}{3 \mu} \Big)\\
        &\leq  2 \exp\Big(- \frac{(\varepsilon N_\ell \bar{p}_k - p_{k,k})^2}{3 N_\ell \bar{p}_k} \Big)\\
        &\leq  2 \exp\Big(- N_\ell \frac{\varepsilon^2 \bar{p}_k}{3} + \frac{2 \varepsilon p_{k,k}}{3} \Big)\\
        &\leq  2 \exp\Big(- N_\ell \frac{\varepsilon^2 \bar{p}}{3} + \frac{2 \varepsilon}{3} \Big)
    \end{align}
    where the second inequality is due to the Chernoff bound (Lemma \ref{lemma:chernoff}).
\end{proof}

Now we can verify the conditions of Theorem \ref{thm:main}.

\begin{proposition}
Let $\rv{G}_{\ell}$ denote the stochastic block model random graph and $\hat{r} := \max_{m\neq n} r_{m,n}$.
We further denote $\bar{p} := \min_{k \in [K]} \bar{p}_k$.
For all $\varepsilon \in (0, \hat{r})$  we have that
\begin{align}
    \forall \ell \in \mathbb{N}:\ \mathbb{P}\Big(\max_{x\in [M]^{N_\ell}} \Delta^{\rv{G}_{N_\ell}}(x) \geq \varepsilon \Big) \leq f_\varepsilon(\ell),
\end{align}
where
\begin{align}
    f_\varepsilon(\ell) := 4 M^{N_\ell + 2} K \exp \Big( - N_\ell^2 \bar{p} \Big(\frac{\varepsilon}{\hat{r}} -  \frac{\varepsilon^2}{\hat{r}^2}\Big)^2 \big/ 12 \Big) + 2 N_\ell \exp\Big( -N_\ell \frac{\varepsilon^2 \bar{p}}{12 \hat{r}} + \frac{\varepsilon}{3 \hat{r}} \Big).
\end{align}
\end{proposition}
\begin{proof}
    We fix any $\ell \in \mathbb{N}$ and denote $\rv{G} := \rv{G}_{\ell}$.
    Inserting the propensity functions for the stochastic block model given in \eqref{eq:block_propensity} and \eqref{eq:propensity_sbm}  yields
    \begin{align}
        \Delta^{\rv{G}}(x) = \max_{(m,k)\to n} r_{m,n} \Big\lvert \frac{1}{N_\ell} \sum_{i: (x_i, s_i) = (m,k)} \frac{d_{i,n}^{\rv{G}}}{d_i^{\rv{G}}} - C_{(m,k)}(x) \frac{\sum_{{k^\prime}\in [K]} C_{(n,{k^\prime})}(x) p_{k,{k^\prime}}}{\bar{p}_k} \Big\lvert.
    \end{align}
    Let $\delta \in (0,1)$ and define the events
    \begin{align}
        \mathcal{A} &:= \Big\{ \max_{x\in [M]^{N}} \Delta^{\rv{G}}(x) \geq \varepsilon \Big\} \\
        \mathcal{B} &:= \Big\{ \forall i: (1-\delta) N_\ell \bar{p}_{s_i} \leq d^{\rv{G}}_i \leq (1+\delta)N_\ell \bar{p}_{s_i}\Big\}.
    \end{align}
    From Lemma~\ref{lemma:sbm_degrees} and the union bound, it follows that
    \begin{align} \label{eq:sbm_sandwich}
        \mathbb{P}(\mathcal{B}^C) \leq 2 N_\ell \exp\Big(-\frac{\delta^2 N_\ell \bar{p}}{3} + \frac{2 \delta}{3}\Big),
    \end{align}
    where $\mathcal{B}^C$ denotes the complement of $\mathcal{B}$.
    We will now show that $\mathbb{P}(\mathcal{A} \cap \mathcal{B})$ is small, from which the Proposition follows when combined with \eqref{eq:sbm_sandwich}.
    We define $\rv{E}^x_{(m,k)\to n}$ as in Lemma \ref{lemma:sbm_counts}.
    For any fixed state $x \in [M]^N$ and any transition $(m,k) \to n$, we have, using the abbreviation $c_{(m,k)} := C_{(m,k)}(x)$,
    \begin{align}
        &\mathbb{P}\Big( r_{m,n} \Big\lvert \frac{1}{N_\ell} \sum_{i : s_{i} = (m,k)} \frac{d^{\rv{G}}_{i,n}(x)}{(1\pm \delta)N_\ell \bar{p}_k} - c_{(m,k)} \frac{\sum_{{k^\prime} \in [K]} c_{(n,k^\prime)} p_{k,{k^\prime}}}{\bar{p}_k}  \Big\rvert \geq \varepsilon \Big)\\
        &= \mathbb{P}\Big( \frac{r_{m,n}}{(1\pm \delta) N_\ell^2 \bar{p}_k}  \big\lvert \rv{E}^x_{(m,k)\to n} - c_{(m,k)} \sum_{{k^\prime} \in [K]} c_{(n,k^\prime)} p_{k,{k^\prime}} (1\pm \delta) N_\ell^2 \big\rvert \geq \varepsilon \Big)\\
        & \leq \mathbb{P}\Big( \big\lvert \rv{E}^x_{(m,k)\to n} - c_{(m,k)} \sum_{{k^\prime} \in [K]} c_{(n,k^\prime)} p_{k,{k^\prime}} N_\ell^2 \big\rvert + \underbrace{c_{(m,k)} \sum_{{k^\prime} \in [K]} c_{(n,k^\prime)} p_{k,{k^\prime}}}_{\leq \bar{p}_k} \delta N_\ell^2\\
        &\qquad \qquad \geq r_{m,n}^{-1}  (1\pm \delta) N_\ell^2 \bar{p}_k \varepsilon \Big)\\
        & \leq \mathbb{P}\Big(\big\lvert \rv{E}^x_{(m,k)\to n} - c_{(m,k)} \sum_{{k^\prime} \in [K]} c_{(n,k^\prime)} p_{k,{k^\prime}} N_\ell^2 \big\rvert \geq r_{m,n}^{-1}  (1\pm \delta) N_\ell^2 \bar{p}_k \varepsilon - \delta N_\ell^2 \bar{p}_k  \Big) \\
        & \leq \mathbb{P}\Big( \big\lvert \rv{E}^x_{(m,k)\to n} - c_{(m,k)} \sum_{{k^\prime} \in [K]} c_{(n,k^\prime)} p_{k,{k^\prime}} N_\ell^2 \big\rvert \geq N_\ell^2 \bar{p} (\hat{r}^{-1} \varepsilon - \hat{r}^{-1} \varepsilon \delta - \delta) \Big).
    \end{align}
    Thus, choosing $\delta = \hat{r}^{-1} \varepsilon / 2$, as we did in the proof of Proposition \ref{prop:ER_f}, adding the maxima, and applying Lemma \ref{lemma:sbm_counts} and the union bound yields
    \begin{align}
        &\mathbb{P}\Big( \max_{x \in [M]^N} \max_{(m,k)\to n} r_{m, n} \Big\lvert \frac{1}{N} \sum_{i : s_{i} = (m,k)} \frac{d^{\rv{G}}_{i,n}(x)}{(1\pm \delta)N_\ell \bar{p}_k} - C_{(m,k)}(x) \frac{\sum_{{k^\prime} \in [K]} C_{(n,{k^\prime})}(x) p_{k,{k^\prime}}}{\bar{p}_k}  \Big\rvert \geq \varepsilon \Big) \nonumber\\
        &\leq 2 M^{N_\ell + 2} K \exp \Big( -\frac{1}{12} N_\ell^2 \bar{p} \Big(\frac{\varepsilon}{\hat{r}} -  \frac{\varepsilon^2}{\hat{r}^2}\Big)^2 \Big)
    \end{align}
    With the same arguments as in the proof of Proposition \ref{prop:ER_f}, this leads to
    \begin{align}
        \mathbb{P}(\mathcal{A}) &\leq \mathbb{P}(\mathcal{A} \cap \mathcal{B}) + \mathbb{P}(\mathcal{B}^C)\\
        &\leq 4 M^{N_\ell + 2} K \exp \Big( -\frac{1}{12} N_\ell^2 \bar{p} \Big(\frac{\varepsilon}{\hat{r}} -  \frac{\varepsilon^2}{\hat{r}^2}\Big)^2 \Big) + 2 N_\ell \exp\Big( -N_\ell \frac{\varepsilon^2 \bar{p}}{12 \hat{r}} + \frac{\varepsilon}{3 \hat{r}} \Big).
    \end{align}
\end{proof}
From the bounding function $f_\varepsilon$ derived above, Theorem~\ref{thm:sbm} follows.

\section{Invariance under graph isomorphism} \label{sec:appendix_isomorph}
Let $\mathcal{G}_N$ denote the set of simple graphs with vertex set~$[N]$.
A graph isomorphism between two simple graphs $G = ([N], E_G)$ and $H = ([N], E_H)$ is a permutation $\tau: [N] \to [N]$ such that
\begin{align}
    (i,j) \in E_G \quad \Leftrightarrow \quad (\tau(i), \tau(j)) \in E_H. 
\end{align}
Hence, we will denote $H = \tau(G)$, if $H$ and $G$ are isomorphic with permutation $\tau$.
Many random graphs that we typically work with are indifferent with respect to the specific node labels, which motivates the following definition:
\begin{definition}
A random graph $\rv{G} \in \mathcal{G}_N$ is called invariant under isomorphism if for any two isomorphic graphs $G, H \in \mathcal{G}_N$ we have
\begin{align}
    \mathbb{P}(\rv{G} = G) = \mathbb{P}(\rv{G} = H).
\end{align}
\end{definition}

\begin{example}
Erd\H{o}s--Rényi random graphs (cf. section \ref{subsec:ER}) are invariant under isomorphism as the probability depends only on the number of edges, which is preserved under graph isomorphism.
Uniformly random $d$-regular graphs (cf. section \ref{subsec:regular}) are also invariant under isomorphism as every $d$-regular graph has equal probability, any not $d$-regular graph has probability $0$, and $d$-regularity is preserved under graph isomorphism.
\end{example}

Let $x \in [M]^N$ be a state and $\tau:[N]\to[N]$ a permutation. We define the permuted state $\tau(x) \in [M]^N$ by $\tau(x)_i := x_{\tau^{-1}(i)}$.
Note that certain observables are identical for $(G,x)$ and $(\tau(G), \tau(x))$, for example the number of edges between nodes of state $m$ and nodes of state $n$.
\begin{definition}
We call a function $f:\mathcal{G}_N \times [M]^N \to \mathbb{R}$ \textit{invariant under isomorphism} if for all permutations $\tau:[N]\to[N]$ and all $(G,x) \in \mathcal{G}_N \times [M]^N$ we have $f(G, x) = f(\tau(G), \tau(x))$.
\end{definition}

\begin{example}
    Let $f(G,x)$ denote the number of edges between nodes of state $m$ and nodes of state $n$, $m\neq n$. Then
    \begin{align}
        f(G,x) &= \sum_{i : x_{i} = m} d^{G}_{i,n}(x)\\
        &= \sum_{i : x_{i} = m} d^{\tau(G)}_{\tau(i),n}(\tau(x))\\
        &= \sum_{i : \tau(x)_{i} = m} d^{\tau(G)}_{i,n}(\tau(x)) = f(\tau(G), \tau(x)),
    \end{align}
    i.e., $f(G,x)$ is invariant under isomorphism.
\end{example}

We consider the following
\begin{proposition}
Let both the random graph $\rv{G} \in \mathcal{G}_N$ and the function $f:\mathcal{G}_N \times [M]^N \to \mathbb{R}$ be invariant under isomorphism.
Let $\tau:[N]\to[N]$ be a permutation and $x \in [M]^N$ a state.
Then we have
\begin{align}
    f(\rv{G}, x) \overset{d}{=} f(\rv{G}, \tau(x)).
\end{align}
\end{proposition}
\begin{proof}
    Define for some fixed $\beta \in \mathbb{R}$
    \begin{align}
        \mathbb{G} &:= \{G \in \set{G}_N\ :\ f(G, x) = \beta\}\\
        \mathbb{G}^* &:= \{G \in \set{G}_N\ :\ f(G, \tau(x)) = \beta\}.
    \end{align}
    Due to the invariance under isomorphism of $f$, we have for any $G \in \mathbb{G}$
    \begin{align}
        \beta = f(G, x) = f(\tau(G), \tau(x))
    \end{align}
    and thus $\tau(G) \in \mathbb{G}^*$.
    Now, let $G^* \in \mathbb{G}^*$. Then
    \begin{align}
        f(\tau^{-1}(G^*), x) = f(G^*, \tau(x)) = \beta
    \end{align}
    and thus $\tau^{-1}(G^*) \in \mathbb{G}$.
    Altogether, we have $\tau(\mathbb{G}) = \mathbb{G}^*$.
    Finally, by the invariance under isomorphism of $\rv{G}$, it follows
    \begin{align}
        \mathbb{P}(f(\rv{G}, x) = \beta) = \mathbb{P}(\rv{G} \in \mathbb{G}) = \mathbb{P}(\rv{G} \in \mathbb{G}^*) = \mathbb{P}(f(\rv{G}, \tau(x)) = \beta).
    \end{align}
\end{proof}

As a consequence, it is sufficient to only deal with ordered system states, for example $x = (1,\dots,1,2,\dots,2,3,\dots) \in [M]^N$, when examining the distribution of $f(\rv{G}, x)$, as any permutation of $x$ would yield an identical distribution. This is exploited in Proposition~\ref{prop:regular_bound}.

\bibliographystyle{elsarticle-num} 
\bibliography{main.bib}

% \end{appendices}
\end{document}